\def\Arxived{arxived}
\newif\ifarxived
\ifx\Arxived\Kaba\else\arxivedtrue\fi
\newif\ifextended
\ifx\Extended\Kaba\else\extendedtrue\fi
\newif\ifprivate
\ifx\Private\Kaba\else\privatetrue\extendedtrue\fi
\newif\ifJapanese
\newif\iftesting
\ifextended\testingtrue\fi
\ifarxived
\extendedtrue
\testingfalse
\privatefalse
\fi
\documentclass[12pt]{article}
\ifextended
\ifarxived
\usepackage[bookmarks=true,backref=page,bookmarksnumbered=true, 
  bookmarkstype=toc,
  pdftitle={Maximality Principle under a Laver-generic supercompact cardinal},%
  pdfauthor={Sakae Fuchino (\quad\quad \quad)},%
  pdfsubject={},%
  pdfkeywords={Generic large cardinal, Laver-generic large cardinal, Maximality Principles}]{hyperref}
\else
\usepackage[dvipdfmx,backref=page,bookmarks=true,bookmarksnumbered=true,
  bookmarkstype=toc,
]{hyperref}
\usepackage[dvipdfmx]{pxjahyper}
\hypersetup{
  colorlinks=true,
  linkcolor=blue,
  filecolor=blue,      
  citecolor=cyan,
  urlcolor=cyan,
  pdftitle={Maximality Principle under a Laver-generic supercompact cardinal},%
  pdfauthor={Sakaé Fuchino (渕野 昌)},%
  pdfsubject={},%
  pdfkeywords={Generic large cardinal, Laver-generic large cardinal, Maximality Principles}
}
\fi 
\else 
\newcommand{\href}[2]{{#2}}
\def\hyperref[#1]#2{{#2}}
\newcommand{\phantomsection}{}
\fi

\ifarxived
\usepackage{graphicx, color}
\else
\usepackage[dvipdfmx]{graphicx, color}
\fi
\usepackage{amsmath, amssymb}
\usepackage{bbm}
\usepackage{dsfont}
\usepackage{bbold}
\usepackage{accents} 
\usepackage{marginnote}
\usepackage[many]{tcolorbox}
\usepackage{mathtools}
\usepackage[normalem]{ulem}
\usepackage{setspace}

\definecolor{darkelectricblue}{rgb}{0.33, 0.39, 0.52}
\definecolor{darkgreen}{rgb}{0.31, 0.47, 0.26}
\ifextended
\newcommand{\extendedcolor}{\color{darkelectricblue}}
\newcommand{\darkgreen}{\color{darkgreen}}

\newcommand{\darkred}{\color[rgb]{0.8,0.1,0.1}}
\newcommand{\darkblue}{\color[rgb]{0.1,0.1,0.8}}
\newcommand{\It}{\it\darkred{}}
\else
\newcommand{\darkgreen}{}
\newcommand{\darkred}{}
\newcommand{\darkblue}{}
\newcommand{\It}{\it{}}
\fi
\usepackage{theorem}
\newcommand{\bbd}[1]{{\mathbb{#1}}}
\theorembodyfont{\ifJapanese\rm\else\it\fi}
\newcount\minute	
\newcount\hour		
\newcount\hourMins  
\ifJapanese
\def\today%
{
  \the\year 年\,\zeroPadTwo{\the\month}月\,\zeroPadTwo{\the\day}日%
}
\fi
\def\now%
{
  \minute=\time    
  \hour=\time \divide \hour by 60 
  \hourMins=\hour \multiply\hourMins by 60
  \advance\minute by -\hourMins 
  \zeroPadTwo{\the\hour}:\zeroPadTwo{\the\minute}%
}
\def\zeroPadTwo#1%
{
  \ifnum #1<10 0\fi    
  #1
}
\title{\scalebox{0.84}[1.6]{\bf Maximality Principle{\ifextended\extendedcolor s and 
      Ressurection Axioms\fi}}\\[2\jot]
\scalebox{0.84}[1.6]{\bf under a Laver-generic {\ifextended{\extendedcolor 
      large{}}\else{supercompact){}}\fi} cardinal}}
\author{\protect\scalebox{1}[1.4]{\ifarxived 
    \ \includegraphics[width=4em]{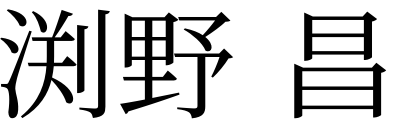}\else 渕野 昌\fi}\medskip\\
  \protect\scalebox{1}[1.4]{\quad Sakaé Fuchino$^{\ast}$}\\
  }
\date{}

\ifarxived
\fi
\ifextended
\fi
\renewcommand{\baselinestretch}{1.2}
\renewcommand{\thefootnote}{(\arabic{footnote})\,}
\iftesting
\fi
\iftesting
\newcommand{\Label}[1]{\label{#1}\marginnote{{\color{cyan}\renewcommand{\baselinestretch}{0.4}\tiny 
		  \rlap{#1}}}}

\else
\newcommand{\Label}[1]{\label{#1}}

\fi

\def\memo#1{\ifprivate\marginnote{{\darkred\normalsize%
      \renewcommand{\baselinestretch}{0.4}\tiny\mbox{}\vspace{-2.52ex}
      \par\relax%
			#1\par\mbox{}}}\else\fi}%
\def\memox#1{}
\def\imemox#1{}

\newcounter{frml}[section]
\newcounter{frmla}[section]
\def\thefrml{{\arabic{section}.\arabic{frml}}}
\def\thefrmla{{$\aleph$\arabic{section}.\arabic{frmla}}}
\def\frmlabel#1{\refstepcounter{frml}{\def\baka{#1}\ifx\baka\empty\else\label{#1}\fi}%
{\rm({\thefrml})\hfill\hfill\hfill}}
\def\ifrmlabel#1{\refstepcounter{frml}{\def\baka{#1}\ifx\baka\empty\else\label{#1}\fi}%
{\iftesting\darkred\fi\rm({\thefrml})\,:\hspace{0.6em}}}
\def\frmlabela#1{\refstepcounter{frmla}{\def\baka{#1}\ifx\baka\empty\else\label{#1}\fi}%
{\rm({\thefrmla})\hfill\hfill\hfill}}
\def\xitem[#1]{\item[\frmlabel{#1}]\mbox{}%
	\iftesting\marginnote{{\renewcommand{%
				\baselinestretch}{0.6}\tiny#1}}\fi\ignorespaces}
\def\xitemq[#1]{\item[\frmlabel{#1}]\mbox{}%
	\ignorespaces}
\def\xitemd[#1]#2{\item[(\ref{#1})$#2$\hfill\hfill\hfill]}
\def\xitemA[#1]{\item[\frmlabela{#1}]\mbox{}%
	\iftesting\marginnote{{\renewcommand{%
				\baselinestretch}{0.6}\tiny#1}}\fi\ignorespaces}
\def\xitemx[#1]{\item[]}
\def\xitemsub[#1]#2{\item[\frmlabel{#1}$_{#2}$]\mbox{}%
	\iftesting\marginnote{{\renewcommand{%
				\baselinestretch}{0.6}\tiny#1}}\fi\ignorespaces}
\def\xxitem[#1][#2]{\item[(\ref{#1}{\makebox[1.4ex][c]{#2}})]\mbox{}%
	\iftesting\marginnote{{\renewcommand{%
				\baselinestretch}{0.6}\tiny\{#1\}\{#2\}}}\fi\ignorespaces}
\def\xitemof#1{{\rm({\ref{#1}})}}
\def\Xitem[#1]{\item[{\makebox[7ex][l]{\rm(\ref{#1})}}]\iftesting\marginnote{{\renewcommand{%
				\baselinestretch}{0.6}\tiny#1}}\fi\ignorespaces}
\def\xitemAof#1{{\rm({\ref{#1}})}}

\newenvironment{xitemize}{\begin{list}{}{\parsep=0.5\smallskipamount%
			\itemindent=-0.4ex%
			\itemsep=0.5\smallskipamount\leftmargin=4em\labelwidth=3em\labelsep=0.7em}}%
							 {\end{list}}
\iftesting
\def\ixitem[#1]{\ifrmlabel{#1}\marginnote{{\color{cyan}\renewcommand{%
				\baselinestretch}{0.6}\qquad\qquad\tiny\rlap{#1}\mbox{}}}\ignorespaces}
\else
\def\ixitem[#1]{\ifrmlabel{#1}\ignorespaces}
\fi
\iftesting
\def\ixitema[#1]{\ifrmlabel{#1}\marginnote{{\color{cyan}\renewcommand{%
				\baselinestretch}{0.6}\tiny\mbox{}\hfill #1}}\ignorespaces}
\else
\def\ixitema[#1]{\ifrmlabel{#1}\ignorespaces}
\fi

\def\assert#1{\noindent\makebox[4.8ex][r]{\rm(\makebox[2.2ex][c]{#1})}\ \ \ignorespaces}
\def\wassert#1{\assert{#1}}
\def\wassertof#1{\makebox[4.8ex][r]{\rm(\makebox[2.2ex][c]{#1})\ }}%
\def\assertof#1{{\rm(#1)}}%

\ifarxived
\usepackage{tikz}
\newcommand*\daimaru[1]{\scalebox{0.6}{\tikz[baseline=(char.base)]{
            \node[shape=circle,draw,inner sep=2pt] (char) {#1};}}}
\else
\def\daimaru#1{\makebox[1em][c]{\mbox{\leavevmode\lower.144ex\hbox{
        \rlap{\hbox to 
          0.76em{\hfil\mbox{}\hfill{}\raisebox{0.054ex}{\scalebox{1.2}{○}}\hfil}}
        \raise0.342ex\hbox to 1em{\hfil{\hspace{0.16em}\footnotesize#1}\hfil}}}}\,}
\fi

\newtheorem{Thm}{\ifJapanese{\bf 定理}\else {\bf Theorem}\fi}[section]
\tcolorboxenvironment{Thm}{
  breakable,
  colback=blue!4!white,
  boxrule=0pt,
  boxsep=0pt,
  left=8pt,right=8pt,top=2pt,bottom=4pt,
  oversize=4pt,
  sharp corners,
  before skip=\topsep,
  after skip=\topsep,
}

\newtheorem{ThmA}{\ifJapanese{\bf 定理\,A\!}\else{\bf Theorem\,A\!}\fi}[section]
\tcolorboxenvironment{ThmA}{
  colback=blue!4!white,
  boxrule=0pt,
  boxsep=0pt,
  left=8pt,right=8pt,top=2pt,bottom=4pt,
  oversize=4pt,
  sharp corners,
  before skip=\topsep,
  after skip=\topsep,
}
\newtheorem{Ex}[Thm]{\ifJapanese{\bf 例}\else {\bf Example}\fi}
\tcolorboxenvironment{Ex}{
  colback=blue!3!white,
  boxrule=0pt,
  boxsep=0pt,
  left=8pt,right=8pt,top=2pt,bottom=4pt,
  oversize=4pt,
  sharp corners,
  before skip=\topsep,
  after skip=\topsep,
}
\newtheorem{Prop}[Thm]{\ifJapanese{\bf 命題}\else{\bf Proposition}\fi}
\tcolorboxenvironment{Prop}{
  colback=blue!4!white,
  boxrule=0pt,
  boxsep=0pt,
  left=8pt,right=8pt,top=2pt,bottom=4pt,
  oversize=4pt,
  sharp corners,
  before skip=\topsep,
  after skip=\topsep,
}

\newtheorem{Lemma}[Thm]{\ifJapanese{\bf 補題}\else{\bf Lemma}\fi}
\tcolorboxenvironment{Lemma}{
  colback=blue!4!white,
  boxrule=0pt,
  boxsep=0pt,
  left=8pt,right=8pt,top=2pt,bottom=4pt,
  oversize=4pt,
  sharp corners,
  before skip=\topsep,
  after skip=\topsep,
}

\tcolorboxenvironment{LemmaA}{
  colback=blue!4!white,
  boxrule=0pt,
  boxsep=0pt,
  left=8pt,right=8pt,top=2pt,bottom=4pt,
  oversize=4pt,
  sharp corners,
  before skip=\topsep,
  after skip=\topsep,
}
\newtheorem{Cor}[Thm]{\ifJapanese{\bf 系}\else{\bf Corollary}\fi}
\tcolorboxenvironment{Cor}{
  colback=blue!4!white,
  boxrule=0pt,
  boxsep=0pt,
  left=8pt,right=8pt,top=2pt,bottom=4pt,
  oversize=4pt,
  sharp corners,
  before skip=\topsep,
  after skip=\topsep,
}

\tcolorboxenvironment{Remark}{
  colback=blue!4!white,
  boxrule=0pt,
  boxsep=0pt,
  left=8pt,right=8pt,top=2pt,bottom=4pt,
  oversize=4pt,
  sharp corners,
  before skip=\topsep,
  after skip=\topsep,
}

\newtheorem{Claim}{{\bf Claim}}[Thm]

\newcommand{\prf}{\noindent\ifJapanese{\bf 証明．\ }\ignorespaces\else{\bf 
		Proof.\ \ }\ignorespaces\fi}

\newcommand{\prfofClaim}{\raisebox{-.4ex}{\Large $\vdash$\ \ }}

\newcommand{\prfof}[1]{\ifJapanese{\bf #1 の証明．\ \ }%
	\ignorespaces\else{\bf Proof of #1:}\ \ \ignorespaces\fi}
\newcommand{\Thmof}[1]{\ifJapanese{定理\,\ref{#1}}\else{Theorem~\ref{#1}}\fi}
\newcommand{\ThmAof}[1]{\ifJapanese{定理\,A\,\ref{#1}}\else{Theorem\,A\,\ref{#1}}\fi}

\newcommand{\bfThmof}[1]{\ifJapanese{\bf 定理\,\ref{#1}}\else{\bf Theorem~\ref{#1}}\fi}
\newcommand{\Lemmaof}[1]{\ifJapanese{補題\,\ref{#1}}\else{Lemma~\ref{#1}}\fi}

\newcommand{\LemmaAof}[1]{\ifJapanese{補題\,A\,\ref{#1}}\else{Lemma\,A\,\ref{#1}}\fi}

\newcommand{\Propof}[1]{\ifJapanese{命題\,\ref{#1}}\else{Proposition~\ref{#1}}\fi}
\newcommand{\bfPropof}[1]{\ifJapanese{\bf 命題\,\ref{#1}}\else{\bf Proposition~\ref{#1}}\fi}
\newcommand{\Corof}[1]{\ifJapanese{系\,\ref{#1}}\else{Corollary~\ref{#1}}\fi}

\newcommand{\Exof}[1]{\ifJapanese{例\,\ref{#1}}\else{Example~\ref{#1}}\fi}
\newcommand{\sectionof}[1]{\ifJapanese{第\ref{#1}節}\else{Section~\ref{#1}}\fi}

\newcommand{\Thmabove}{{\ifJapanese 定理\else Theorem\fi\ \number\theThm}}

\newcommand{\ubecause}[3]{\underbrace{{}#1{}%
  \ifx\bakakaba#2\bakakaba\rule[-0.72ex]{0pt}{1pt}\else\rule[#2]{0pt}{1pt}\fi}_{\mbox{\footnotesize\clap{#3}}}}
\newcommand{\obecause}[3]{\overbrace{{}#1{}%
  \ifx\bakakaba#2\bakakaba\rule[1.62ex]{0pt}{1pt}\else\rule[#2]{0pt}{1pt}\fi}^{\mbox{\footnotesize\clap{#3}}}}
\newsavebox{\qedbox}\sbox{\qedbox}{
{\unitlength=0.05mm \begin{picture}(40,60)
\put(0,0){\framebox(30,44)[cc]{}}
\put(30,-7){\rule{7\unitlength}{44\unitlength}}
\put(10,-7){\rule{27\unitlength}{7\unitlength}}
\end{picture}}}
\newcommand{\qed}{\mbox{}\hfill\usebox{\qedbox}}
\newcommand{\smallqed}%
{\mbox{}\smallskip\hfill\raisebox{-.4ex}{\Large $\dashv$}}
\newcommand{\qedof}[1]%
{\mbox{} \hspace*{\fill}{\usebox{\qedbox}{\tiny~(#1)}}}
\newcommand{\Qedof}[1]%
{\mbox{} \hspace*{\fill}{\usebox{\qedbox}%
{\tiny~(#1~\number\theThm)}}}
\newcommand{\QedAof}[1]%
{\mbox{} \hspace*{\fill}{\usebox{\qedbox}%
{\tiny~(#1~\number\theThmA)}}}
\newcommand{\qedofThm}{\Qedof{\ifJapanese 定理\else Theorem\fi}}

\newcommand{\qedofCor}{\Qedof{\ifJapanese 系\else Corollary\fi}}
\newcommand{\qedofProp}{\Qedof{\ifJapanese 命題\else Proposition\fi}}
\newcommand{\qedofLemma}{\Qedof{\ifJapanese 補題\else Lemma\fi}}

\newcommand{\qedskip}{\medskip}
\newcommand{\qedofEx}{\Qedof{\ifJapanese 例\else Example\fi}}

\newcommand{\qedofClaim}%
{\mbox{}\hfill\raisebox{-.4ex}{\Large $\dashv$ }\nolinebreak%
\mbox{\tiny~(Claim~\number\theClaim)}}
\newcommand{\qedofClaimA}%
{\mbox{}\hfill\raisebox{-.4ex}{\Large $\dashv$ }\nolinebreak%
\mbox{\tiny~(Claim~A\,\number\theClaimA)}}
\newcommand{\qedofClaimAof}[1]%
{\mbox{}\hfill\raisebox{-.4ex}{\Large $\dashv$ }\nolinebreak%
\mbox{\tiny~(Claim~A\,\ref{#1})}}
\newcommand{\qedofSubclaim}%
{\mbox{}\hfill\raisebox{-.4ex}{\Large $\dashv$ }\nolinebreak%
\mbox{\tiny~(Subclaim~\number\theSubclaim)}}
\newcommand{\cardof}[1]{\mathopen{|\,}#1\mathclose{\,|}}

\newcommand{\Card}{{\it Card\/}}

\newcommand{\setof}[2]{\{#1\,:\,#2\}}
\newcommand{\ssetof}[1]{\{#1\}}

\newcommand{\subseteqand}[1]{\mathrel{\mathop{\subseteq}%
		\limits_{\scriptscriptstyle\hbox to 14pt{$\scriptscriptstyle #1$\hss}}}}

\newcommand{\mapping}[3]{#1:#2\rightarrow #3}

\newcommand{\elembed}[3]{#1:#2\stackrel{\preccurlyeq\hspace{0.8ex}}{\rightarrow}#3}
\newcommand{\Elembed}[4]{#1:#2\stackrel{\prec\hspace{0.8ex}}{\rightarrow}_{#4}#3}

\newcommand{\fnsp}[2]{\mbox{}^{{#1}\hspace{-0.02em}}#2}
\newcommand{\imageof}{{}^{\,{\prime}{\prime}}}

\newcommand{\seqof}[2]{\langle#1\,:\,#2\rangle}
\newcommand{\pairof}[1]{\langle#1\rangle}
\newcommand{\psof}[1]{{\mathcal P}\/(#1)}
\newcommand{\forces}[2]{\,\|\hspace{-.35ex}\mbox{\sf--}_{\,#1\,}%
\mbox{\rm``}\,#2\,\mbox{\rm''}}
\newcommand{\notforces}[2]{\rlap{\rm\ 
 /}\|\hspace{-.35ex}\mbox{\sf--}_{\,#1\,}%
 \mbox{\rm``}\,#2\,\mbox{\rm''}}
\newcommand{\decides}[2]{\,\|_{\,#1\,}%
\mbox{\rm``}\,#2\,\mbox{\rm''}}

\newcommand{\modelof}[1]{\models\!\mbox{\rm``\,}#1\mbox{\rm''}}

\newcommand{\crit}{\mbox{\it crit\/}}

\newcommand{\bbone}{{\mathord{\mathbb{1}}}}

\newcommand{\circleq}{\mathrel{{\leqslant}%
		\hspace{-0.86ex}{\lower-0.53ex\hbox{$\scriptscriptstyle\circ$}}}}

\newcommand{\restr}{\restriction}

\newcommand{\Col}{{\rm Col}}

\newcommand{\Fn}{{\rm Fn}}

\newcommand{\trcl}{\mathop{\mbox{\it trcl\/}}}

\newcommand{\natnums}{{\bbd{N}}}

\newcommand{\poP}{\bbd{P}}
\newcommand{\poQ}{\bbd{Q}}
\newcommand{\poR}{\bbd{R}}
\newcommand{\On}{{\rm On}}

\newcommand{\genG}{\mathbb{G}}

\newcommand{\utpoP}{\utilde{\mathbb{P}}}
\newcommand{\utildepoQ}{\utilde{\mathbb{Q}}}
\newcommand{\utpoQ}{\utilde{\mathbb{Q}}}
\newcommand{\utpoR}{\utilde{\mathbb{R}}}

\newcommand{\genH}{\mathbb{H}}

\newcommand{\condp}{\mathbbm{p}}

\newcommand{\LT}{{<}\,}
\newcommand{\LE}{{\leq}\,}
\newcommand{\GT}{{>}\,}

\newcommand{\ctenten}{,\mbox{}\hspace{0.08ex}{.}{.}{.}\hspace{0.1ex}}
\newcommand{\ctentenc}{,{}\linebreak[0]\hspace{0.04ex}{{.}{.}{.}\hspace{0.1ex},\,}\linebreak[0]}

\newcommand{\xmbox}[1]{ $\relax{\rm #1}\relax$ }

\newcommand{\continuum}{2^{\aleph_0}}

\newcommand{\calC}{{\mathcal C}}
\newcommand{\calD}{{\mathcal D}}

\newcommand{\calH}{{\mathcal H}}

\newcommand{\calL}{{\mathcal L}}
\newcommand{\calM}{{\mathcal M}}

\newcommand{\calP}{{\mathcal P}}

\newcommand{\calS}{{\mathcal S}}

\newcommand{\calU}{{\mathcal U}}

\newcommand{\uta}{\utilde{a}}
\newcommand{\utb}{\utilde{b}}

\newcommand{\Lin}{{\calL}_{\in}}

\newcommand{\ZFC}{{\sf ZFC}}
\newcommand{\CH}{{\sf CH}}

\newcommand{\SCH}{{\sf SCH}}

\newcommand{\MA}{{\sf MA}}
\newcommand{\MM}{{\sf MM}}
\newcommand{\MP}{{\sf MP}}
\newcommand{\RA}{{\sf RA}}
\newcommand{\UR}{{\sf UR}}
\newcommand{\TUR}{{\sf TUR}}
\newcommand{\BfRA}{{\mathbb{R}\mathbb{A}}}
\newcommand{\PFA}{{\sf PFA}}

\newcommand{\FRP}{{\sf FRP}}

\newcommand{\GRP}{{\sf GRP}}
\newcommand{\RC}{{\sf RC}}

\newcommand{\refl}{{\mathfrak{r}\mathfrak{e}\mathfrak{f}\mathfrak{l}\,}}

\newcommand{\SDLS}{{\sf SDLS}}
\newcommand{\II}{{\mathrm{II}}}

\newcommand{\st}{such that}
\newcommand{\wrt}{with respect to}
\newcommand{\Wolog}{Without loss of generality}

\newcommand{\Tfae}{The following are equivalent}
\newcommand{\tfae}{the following are equivalent}

\newcommand{\po}{poset}
\newcommand{\pos}{posets}
\newcommand{\uniV}{{\sf V}}
\newcommand{\uniL}{{\sf L}}
\newcommand{\uniW}{{\sf W}}

\newcommand{\Pkl}[2]{\ifx\bakakaba#1\bakakaba\ifx\bakakaba#2\bakakaba{\mathcal 
    P}_\kappa(\lambda)\else{\mathcal P}_\kappa(#2)\fi\else{\mathcal P}_{#1}(#2)\fi}

\newcommand{\utildeT}[1]{%
  \underaccent{{\sim}}{#1}}
\newcommand{\utildeS}[1]{%
	\hbox to 0pt{\smash{$\mathop{\scriptstyle #1}\limits_{%
				\raisebox{0.6ex}[0pt]{$\scriptscriptstyle\sim$}}$}\hss}%
	\relax\phantom{\mathord{{#1}_{\rule[-0.6ex]{0pt}{1pt}}}}}
\newcommand{\utildeSS}[1]{%
	\hbox to 0pt{$\mathop{\scriptscriptstyle #1}%
		\limits_{\scriptscriptstyle\sim}$\hss}%
		\relax\phantom{\underline{#1}}}
\newcommand{\utilde}[1]{%
	\mathchoice{\utildeT{#1}}{\utildeT{#1}}{\utildeS{#1}}{\utildeSS{#1}}}

{\end{minipage}\end{trivlist}}

\begin{document}
\maketitle
\renewcommand{\thefootnote}{$\ast$\ }
  \footnotetext{Graduate School of System Informatics, Kobe University \\Rokko-dai 1-1, Nada, Kobe 657-8501 Japan
   \\
    \quad\scalebox{0.95}[1]{\tt fuchino@diamond.kobe-u.ac.jp}}

\ifextended
\phantomsection
\addcontentsline{toc}{section}{* Maximality Principle under a 
  Laver-generic supercompact cardinal} 
\addcontentsline{toc}{section}{**** by S.Fuchino}
\fi

\ifextended
\addcontentsline{toc}{section}{Abstract}
\fi
\begin{abstract}
\ifextended{
  \extendedcolor
  Set-theoretic axioms formulated in terms of existence of a Laver-generic large cardinal 
  were introduced in \cite{sfetal-II} and studied further in \cite{sfetal-III}, \cite{FuOt}, 
  \cite{fuchino-sakai}.   
  These axioms, let us call them Laver-genericity axioms,  claim the existence 
  of a $\calP$-Laver generic large cardinal for various classes $\calP$ of proper or 
  semi-proper \pos, and they still vary depending on the notions of large cardinal 
  involved, and a modification (tightness) of the definition of Laver-genericity.  

  Laver-genericity axioms we consider here are divided into three groups depending on 
  whether they imply that 
  the Laver generic large cardinal $\kappa$ is
  $\aleph_2=(2^{\aleph_0})^+$, or it is $\aleph_2 = 2^{\aleph_0}$, or else it is very large and
  $= 2^{\aleph_0}$ (see the Trichotomy Theorem (\Thmof{p-laver-1})).

	Many set-theoretic axioms and principles considered in the recent development of set theory 
	follow from a Laver-genericity axiom in one of these three groups, and by this, 
  they are placed uniformly in a global context (see \hyperref[figure3]{Figure 3}).  

  In spite of this very strong unifying feature of the Laver genericity axioms, we show 
  that Maximality  
  Principle (\MP) without 
  parameters is independent over \ZFC\ with any 
  of the Laver-genericity axioms we consider in our present context (\Thmof{p-max-5}, 
  \Thmof{p-indep-0-7}). Similar 
  independence is also shown for parameterized versions of Maximality Principles
  (\Thmof{p-para-max-0}, \Thmof{p-para-max-2}). 

  In contrast to these independence results, we can show that local versions of Maximality 
  Principle as well as versions of Resurrection Axioms including the Unbounded Resurrection 
  Axioms of Tsaprounis follow from the 
  existence of a tightly Laver-generic large cardinal for a strong enough notion of large 
  cardinal (\Thmof{p-UR-4}, \Thmof{T-Lg-RA-0}, \Thmof{p-resurr-0}). 
}\else
We give a survey on the set-theoretic axioms formulated in terms of existence of a 
Laver-generic large cardinal.

We show that the Maximality Principle without parameters is independent over \ZFC\ with the axiom 
asserting the existence of a $\calP$-Laver generically supercompact cardinal for an 
iterable class of \pos\ $\calP$ as far as the existence of such a cardinal can be forced  
naturally starting from a genuine supercompact cardinal.  
\fi
\end{abstract}

\ifextended
{\extendedcolor 
\phantomsection
\addcontentsline{toc}{section}{Contents}
\newcommand{\myscalebox}[1]{\scalebox{0.88}[1.06]{#1}}
\begin{quotation}
	\footnotesize
	\noindent
	\centerline{
      \normalsize\tt\quad\ Contents\hspace{6em}\mbox{}}\mbox{}\\
  {\mbox{}\hspace{-1.6em}\tt\makebox[3.4ex][l]{\ref{intro}.}%
   \hyperref[intro]{\tt\myscalebox{Introduction}}}\ \ \dotfill\ \ {\pageref{intro}}\\ 
  {\mbox{}\hspace{-1.6em}\tt\makebox[3.4ex][l]{\ref{gen}.}%
   \hyperref[gen]{\tt\myscalebox{Generic large cardinals}}}\ \ \dotfill\ \ {\pageref{gen}}\\ 
  {\mbox{}\hspace{-1.6em}\tt\makebox[3.4ex][l]{\ref{laver}.}%
    \hyperref[laver]{\tt\myscalebox{Laver-generic large cardinals}}}\ \ \dotfill\ \ {\pageref{laver}}\\ 
  {\mbox{}\hspace{-1.6em}\tt\makebox[3.4ex][l]{\ref{max}.}%
    \hyperref[max]{\tt\myscalebox{Maximality Principle}}}\ \ 
  \dotfill\ \ {\pageref{max}}\\   
  {\mbox{}\hspace{-1.6em}\tt\makebox[3.4ex][l]{\ref{indep}.}%
    \hyperref[indep]{\tt\myscalebox{Independence of \MP\ under a Laver-gen.\ 
        large cardinal}}}\ \ \dotfill\ \ {\pageref{indep}}\\    
  {\mbox{}\hspace{-1.6em}\tt\makebox[3.4ex][l]{\ref{para-max}.}%
    \hyperref[para-max]{\tt\myscalebox{Boldface Maximality Principle for an iterable class $\calP$ of \pos}}}\\
  {\phantom{\mbox{}\hspace{-1.6em}\tt\makebox[3.4ex][l]{??.}}%
    \hyperref[para-max]{\tt\myscalebox{\quad and Laver-genericity}}}\ \ \dotfill\ \ {\pageref{para-max}}\\   {\mbox{}\hspace{-1.6em}\tt\makebox[3.4ex][l]{\ref{resurr}.}%
    \hyperref[resurr]{\tt\myscalebox{Resurrection Axioms}}}\ \ \dotfill\ \ {\pageref{resurr}}\\    
  {\mbox{}\hspace{-1.6em}\hyperref[ref]{\tt References}}\ \ \dotfill\ \ 
  {\pageref{ref}}\medskip\\ 
  {\mbox{}\hspace{0em}\hyperref[figure1]{\tt Figure 1}}\ \ \dotfill\ \ 
  {\pageref{figure1}}\\ 
  {\mbox{}\hspace{0em}\hyperref[figure2]{\tt Figure 2}}\ \ \dotfill\ \ 
  {\pageref{figure2}}\\ 
  {\mbox{}\hspace{0em}\hyperref[figure3]{\tt Figure 3}}\ \ \dotfill\ \ 
  {\pageref{figure3}}\\ 

\end{quotation}}
\fi
\renewcommand{\thefootnote}{}
\footnotetext{{\it Date:} May 8, 2023
  \qquad {\it Last update:} 
  \today\ (\now\ JST)\vspace{-1\smallskipamount}
}
\footnotetext{{\it MSC2020 Mathematical Subject Classification:} 03E35, 03E50, 03E55, 03E37
  	\vspace{-1\smallskipamount}}
\footnotetext{\mbox{}\\[-2.52ex]The research is supported by Kakenhi Grant-in-Aid for 
  Scientific Research (C) 20K03717. \\The author would like to thank Kaethe Minden for making 
  him aware of \cite{minden}. He also would like to thank Gunter Fuchs, Takehiko Gappo, 
  Paul Larson, Hiroshi Sakai, and Kostas Tsaprounis for many helpful remarks and comments.} 

\ifextended
\ifprivate
\footnotetext{This is a private extended version of the paper with the 
  the title ``Maximality Principles under a Laver-generic supercompact cardinal'' to appear 
  in RIMS Kôkyûroku No.????. \memo{fill the ???? in all versions}
  \par All additional 
  details not contained in the submitted version of the paper are either typeset in 
  {\extendedcolor ``dark electric blue''} in case the text is also included in the 
  (non-private) extended 
  version of the paper, or in {\darkred ``dark red''} if the comment are thought only for the 
  private version. 
  \sout{The numbering of the assertions is kept identical with the submitted version.}
  Since the changes from the submitted version are now quite extensive, I am not trying to keep the 
  numbering of the theorems and assertions identical with the numbering in the submitted version. 

  The most up-to-date pdf-file of this private extended version is downloadable as:\\
  \ifprivate\href{https://fuchino.ddo.jp/papers/RIMS2022-RA-MP-xx.pdf}{\tt 
  https://fuchino.ddo.jp/papers/RIMS2022-RA-MP-xx.pdf}\else
  \href{https://fuchino.ddo.jp/papers/RIMS2022-RA-MP-x.pdf}{\tt 
  https://fuchino.ddo.jp/papers/RIMS2022-RA-MP-x.pdf}\fi 

  Some of the materials in this extended version will be used in the forthcoming 
  \cite{future}. 
}
\else
\footnotetext{\extendedcolor This is an extended version of the paper with the 
  title ``Maximality Principles under a Laver-generic supercompact cardinal'' to appear in 
  RIMS Kôkyûroku.  
  \par\extendedcolor  All additional 
  details not contained in the submitted version of the paper are either typeset in 
  dark electric blue (the color in which this paragraph is typeset) or put in a separate appendices. 
 \sout{The numbering of the assertions is kept identical with the submitted version.}
  Since the changes from the submitted version are now quite extensive, I am not trying to keep the 
  numbering of the theorems and assertions identical with the numbering in the submitted version. 

  The most up-to-date pdf-file of this extended version is downloadable as:\medskip\\
  \qquad\qquad\href{https://fuchino.ddo.jp/papers/RIMS2022-RA-MP-x.pdf}{\tt 
  https://fuchino.ddo.jp/papers/RIMS2022-RA-MP-x.pdf} \medskip\\
  The materials in this extended version may be reused in the forthcoming \cite{future}.
}
\fi
\else
\footnotetext{A pdf file of updated and extended version of this paper possibly with more 
  details and 
  proofs is downloadable as:\qquad \href{https://fuchino.ddo.jp/papers/RIMS2022-RA-MP-x.pdf}{\tt https://fuchino.ddo.jp/papers/RIMS2022-RA-MP-x.pdf}}
\fi

\renewcommand{\thefootnote}{\arabic{footnote})\,}

\section{Introduction}
\Label{intro}
{\ifextended\extendedcolor
Set-theoretic axioms formulated in terms of existence of a Laver-generic large cardinal 
(see \sectionof{laver} for definition) 
were introduced in \cite{sfetal-II} and studied further in \cite{sfetal-III}, \cite{FuOt}, 
\cite{fuchino-sakai}.   

More precisely, these axioms --- let us call them here {\It Laver-genericity axioms} --- claim the 
existence  
of a $\calP$-Laver generic large cardinal for various classes $\calP$ of proper or 
semi-proper \pos, and they still vary depending on the notions of large cardinal involved, 
and a modification (tightness) of Laver-genericity.  

We restrict ourselves here to classes $\calP$ of \pos\ which are proper or semi-proper since we 
want to have axioms which imply (or at least compatible with) various reflection principles, see 
\sectionof{gen} and \sectionof{laver}. 

Laver-genericity axioms are divided into three groups depending on whether they imply that 
the Laver generic large cardinal $\kappa$ whose existence is claimed by the axioms is
$\aleph_2=(2^{\aleph_0})^+$, or it is $\aleph_2 = 2^{\aleph_0}$, or else it is very large and
$= 2^{\aleph_0}$ (see the Trichotomy Theorem (\Thmof{p-laver-1})).

By this trichotomy, the Laver-generic large cardinal we 
consider here is proved to be unique under the respective Laver-genericity axiom. 

Many set-theoretic axioms and principles considered in the recent development of set theory 
follow from a Laver-genericity axiom in one of these three groups, and by this, 
they are placed uniformly in a global context (see \hyperref[figure3]{Figure 3}).  

In sections \ref{gen}, \ref{laver} of 
the present note, we give an improved and streamlined presentation of the Laver-genericity 
axioms. Most of the 
materials presented in these sections are already stated in \cite{sfetal-I} or 
\cite{sfetal-II} but there are also a couple of improvements and new results.
The extended version of the paper {\ifextended\extendedcolor you are reading now{}\fi} also contains detailed proofs of the results mentioned in these sections. 

In spite of the very strong unifying feature of the Laver genericity axioms, we can show that 
Maximality   
Principle (\MP) without 
parameters is independent over \ZFC\ with any 
of the Laver-genericity axioms we consider in our present context (\Thmof{p-max-5}, 
\Thmof{p-indep-0-7}). Similar 
independence is also shown for parameterized versions of Maximality Principles
(\Thmof{p-para-max-0}, \Thmof{p-para-max-2}). 

In contrast to these independence results, the local versions of Maximality 
Principle as well as versions of Resurrection Axioms including the Unbounded Resurrection 
Axioms of Tsaprounis are consequences of the 
Laver-genericity axioms for a strong enough notion of large 
cardinal (\Thmof{p-UR-4}, \Thmof{T-Lg-RA-0}, \Thmof{p-resurr-0}). 

In the following we are working in the framework of \ZFC. All classes are definable by some
$\Lin$-formula where $\Lin$ is the language of \ZFC\ consisting solely of the $\in$-symbol. 
Sometimes the language is extended with  a constant symbol or a unary relation symbol. This 
is in particular the case when we are talking about ``$V_\delta\prec\uniV$'' or that ``there are 
stationarily many $\delta$ with certain large cardinal property and 
\st\ $V_\delta\prec\uniV$''. Even in such cases \ZFC\ is meant the axiom system in the 
original language $\Lin$. 

Regardless of this convention, we sometimes choose a narration which may sound that 
we would be  
working in some higher order set theory. This happens in particular when we are talking 
about the notions of Laver-generic large cardinal: we may do this since it is proved in 
\cite{fuchino-sakai} 
that the notions of Laver generic cardinals are actually formalizable in the context of \ZFC.  
\fi}
\section{Generic large cardinals}
\Label{gen}
Let us begin with recalling the definition of supercompact cardinal: A cardinal $\kappa$ 
is {\It supercompact}\/ if, for any $\lambda>\kappa$, there are classes $j$, $M$ \st\quad
\ifarxived\daimaru{1} \else ① \fi $\Elembed{j}{\uniV}{M}{\kappa}$,
\quad\ifarxived\daimaru{2} \else ② \fi $j(\kappa)>\lambda$ and
\quad\ifarxived\daimaru{3} \else ③ \fi$\fnsp{\lambda}{M}\subseteq M$.  

Here, ``$\darkred\Elembed{j}{N}{M}{\kappa}$'' denotes the set of conditions that $N$ 
and $M$ are transitive (sets or classes); $j$ is a non-trivial elementary embedding of  the 
structure $(N,\in)$ into the structure $(M,\in)$; $\kappa\in N$, and
$\crit(j)=\kappa$.

Note that a supercompact cardinal is a large large cardinal which is a normal measure one 
limit of measurable cardinals{\ifextended\extendedcolor\ (see e.g. \cite{higher-inf} 
  Proposition 22.1)\fi}, and more. This is not the case with the generic large 
cardinal version of the notion of supercompactness (e.g. see Examples \ref{ex-gen-1}, 
\ref{ex-gen-2} below).

For a class $\calP$ of \pos, a cardinal $\kappa$ 
is {\It$\calP$-generically supercompact}\/ ({\It$\calP$-gen.\ supercompact}, for short) if, 
for  
every $\lambda>\kappa$,  there is 
$\poP\in\calP$ \st, for $(\uniV,\poP)$-generic filter $\genG$, there are $j$,
$M\subseteq\uniV[\genG]$ \st\quad\ifarxived\daimaru{1} \else ① \fi
$\Elembed{j}{\uniV}{M}{\kappa}$,\quad\ifarxived\daimaru{2} \else ②  \fi
$j(\kappa)>\lambda$, 
and\quad\ifarxived\daimaru{3}' \else ③' \fi$j\imageof{\lambda}\in M$.

{\ifextended\extendedcolor
In case of genuine supercompactness, the condition $j\imageof\lambda\in M$ is equivalent to
$\fnsp{\lambda}{M}\subseteq M$ for $M$ obtained as the ultrapower of $V$ by 
an $\omega_1$-complete ultrafilter (see Kanamori \cite{higher-inf}, Proposition 22.4). In 
general we do not have this equivalence for generic supercompactness. However this condition still 
implies certain closedness of $M$: 

\begin{Lemma}\extendedcolor{\rm (Lemma 2.5 in \cite{sfetal-II})}
  \Label{L-lt-conti-0}
  Suppose that $\genG$ is a $(\uniV,\poP)$-generic filter for a
  \po\/ $\poP\in\uniV$, and
  $\Elembed{j}{\uniV}{M\subseteq\uniV[\genG]}{\kappa}$ for a cardinal $\kappa$ is \st, for  
  a cardinal in $\uniV$ $\lambda$ with $\kappa\leq\lambda$, we have $j\imageof\lambda\in M$. 

  \assert{1} For any set $A\in\uniV$  
  with $\uniV\models\cardof{A}\leq\lambda$, we have $j\imageof A\in M$. \smallskip

  \assert{2} $j\restr\lambda$, $j\restr\lambda^2\in M$.\smallskip

  \assert{3} For any $A\in\uniV$ with $A\subseteq\lambda$ or $A\subseteq\lambda^2$ 
  we have $A\in M$.\smallskip

  \assert{4} $(\lambda^+)^M\geq(\lambda^+)^\uniV$, Thus, if
  $(\lambda^+)^\uniV=(\lambda^+)^{\uniV[\genG]}$,  then
  $(\lambda^+)^M=(\lambda^+)^\uniV$. \smallskip

  \assert{5} $\calH(\lambda^+)^\uniV\subseteq M$.\smallskip

  \assert{6} $j\restr A\in M$ for all $A\in\calH(\lambda^+)^\uniV$. 
\end{Lemma}
\prf \assertof{1}: In $\uniV$, let $\mapping{f}{\lambda}{A}$ be a surjection. 

For each $a\in A$ with $a=f(\alpha)$, we have
\begin{xitemize}
\xitem[pr-L-lt-conti-0] 
  $j(a)=j(f(\alpha))=j(f)(j(\alpha))$
\end{xitemize}
by elementarity. 
Thus $j\imageof A=j(f)\imageof(j\imageof\lambda)$. Since $j(f)$,
$j\imageof\lambda\in M$, it follows that $j\imageof A\in M$.\smallskip

\assertof{2}: Since $j\imageof\lambda\in M$ and $(j\restr\lambda)(\xi)$ for 
$\xi\in\lambda$ is the $\xi$th element of $j\imageof\lambda$, $j\restr\lambda$ is 
definable subset of $\lambda\times j\imageof\lambda$ in $M$ and hence is an 
element of $M$. Similarly, $j\restr\lambda^2\in M$. \smallskip

\assertof{3}: Suppose that $A\in\uniV$ and $A\subseteq\lambda$ (the case of
$A\subseteq\lambda^2$ can be treated similarly). Then $j\imageof A\in M$ by 
\assertof{1}. Thus, by \assertof{2},
$A=(j\restr\lambda)^{-1}\imageof(j\imageof A)\in M$. \smallskip

\assertof{4}: Suppose that $\mu<(\lambda^+)^\uniV$. Then there is $A\in\uniV$ 
with $A\subseteq\lambda^2$ \st\ $A$ codes the order type of $\mu$. $A\in M$ by 
\assertof{3}. Thus $M\modelof{\cardof{\mu}\leq\lambda}$.

If $(\lambda^+)^\uniV=(\lambda^+)^{\uniV[\genG]}$, we have
\begin{xitemize}
\xitem[pr-L-lt-conti-1] 
  $(\lambda^+)^\uniV=(\lambda^+)^{\uniV[\genG]}\geq(\lambda^+)^M\geq(\lambda^+)^\uniV$.
\end{xitemize}

\assertof{5}: For $A\in\calH(\lambda^+)^\uniV$, let $U\in\uniV$ be \st\
$\trcl(A)\subseteq U$ and $\uniV\modelof{\cardof{U}=\lambda}$. Let 
$c_A\subseteq\lambda^2$ and $d_A$, $e_A\subseteq\lambda$ be \st\ $c_A$, $d_A$, $e_A\in\uniV$ 
and 
\begin{xitemize}
\xitem[pr-L-lt-conti-2] 
  $\pairof{\lambda, c_A, d_A, e_A}\cong\pairof{U,\in\restr U^2, \trcl(A), A}$. 
\end{xitemize}
By \assertof{3}, $c_A$, $d_A$, $e_A\in M$ and hence 
$\pairof{\lambda, c_A, d_A, e_A}\in M$. Since $\trcl(A)$ and 
then $A$ can be recovered 
from this quadruplet in $M$, it follows that $A\in M$. 
\smallskip

\assertof{6}: Suppose that $A\in\calH(\lambda^+)^\uniV$. Since $A\in M$ by \assertof{5}, 
it is enough to show 
that $j\restr\trcl(A)\in M$. 

We have $\trcl(A)\in\calH(\lambda^+)^\uniV$ and hence $\trcl(A)\in M$ by 
\assertof{5}. Thus $j\imageof\trcl(A)$, $j\imageof(\in\restr\trcl(A))\in M$ by \assertof{1}. 
But then the mapping $(j\restr\trcl(A))^{-1}$ is the Mostowski collapse of $j\imageof\trcl(A)$. 
Thus $j\restr\trcl(A)\in M$. 
\qedofLemma
\qedskip
\fi}

\begin{Ex}
  \Label{ex-gen-1}
  Suppose that $\kappa$ is a supercompact cardinal and
  $\calP=\Col(\aleph_1,\kappa)$ (the standard collapsing of all cardinals strictly between
  $\aleph_1$ and $\kappa$ by countable conditions). 
  Then for a $(\uniV,\poP)$-generic $\genG$, we have 
  $\kappa=(\aleph_2)^{\uniV[\genG]}$ and
  $\uniV[\genG]\modelof{\kappa\mbox{ is }\sigma\xmbox{-closed-gen.\ super\-compact}}$. \qed
\end{Ex}
\begin{Ex}
  \Label{ex-gen-2}
  If \MA\ is forced starting from an 
  supercompact cardinal $\kappa$ with an ccc-iteration of length $\kappa$ in finite support 
  along with a 
  supercompact Laver function, then we obtain a model in which $\kappa$ is the continuum
  (though still quite large, e.g. hyper-hyper etc.\ weakly Mahlo, and more) and it is 
  ccc-gen.\ supercompact in the generic extension. \qed
\end{Ex}

These examples will be revisited in \Thmof{p-laver-0} below. The situation created in 
\Exof{ex-gen-1} can be also seen as a strong reflection property. 

\begin{Thm}\Label{p-gen-0}{\rm(B.\,K\"onig \cite{koenig}\,)} \Tfae:\medskip

  \wassert{a} Game Reflection 
  Principle (\GRP) holds.\smallskip

  \wassert{b} $\aleph_2$ is $\sigma$-closed-gen.\ supercompact.\qed
\end{Thm}

As in \cite{sfetal-I},  what we call the {\darkred Game Reflection Principle} 
({\darkred\GRP}) is the principle called $\GRP^+$ in \cite{koenig}. As its name suggests, 
\GRP\ is actually a  
reflection statement about the non-existence of winning strategy of certain 
games of length $\omega_1$ down to subgames of size $\LT\aleph_2$.

We will not go into the details of the definition of \GRP\ but just note that \GRP\ 
implies the Continuum Hypothesis (\CH) and it implies practically all reflection principles with 
reflection down to $\LT\aleph_2$ available under \CH:

\begin{xitemize}
\xitem[x-gen-0] \GRP\ implies Rado's Conjecture (\RC) (K\"onig, \cite{koenig}). 

\xitem[x-gen-1] \GRP\ implies strong downward L\"owenheim-Skolem Theorem of
  $\calL^{\aleph_0,\II}_{stat}$ down to $\LT\aleph_2$
  ($\SDLS(\calL^{\aleph_0,\II}_{stat},\LT\aleph_2)$  in the notation of 
  \cite{sfetal-I}). 

\xitem[x-gen-2] 
  Each of \RC\ and $\SDLS^-(\calL^{\aleph_0}_{stat},\LT\aleph_2)$ (a weakening of
  $\SDLS(\calL^{\aleph_0,\II}_{stat},\LT\aleph_2)$) 
  implies Fodor-type 
  Reflection Principle (\FRP), see \cite{sf-note} and \cite{sfetal-I}. 
\xitem[x-gen-2-0] 
  \RC\ implies a strong form of Chang's Conjecture (Todor\v{c}evi\'{c}, \cite{stevo1993})
\xitem[x-gen-3] 
\FRP\ is known to be equivalent to many ``mathematical'' 
  reflection principles with reflection down to $<\aleph_2$, see \cite{fuchino-left}, 
  \cite{balogh}, \cite{pre-hilbert}, \cite{fjetal}, \cite{more}. 
\xitem[x-gen-4] \GRP\ implies a "generic" solution to the Hamburger's problem 
  (see \Corof{p-gen-2} below, for the original Hamburger's Problem see \cite{sfetal-III} 
  and reference 
  given there).  
\end{xitemize}

Some of these and some other implications are put together in the following diagram:\bigskip
\memo{!!! references in the diagram should be renumbered.}
\newpage

\includegraphics[width=0.94\textwidth]{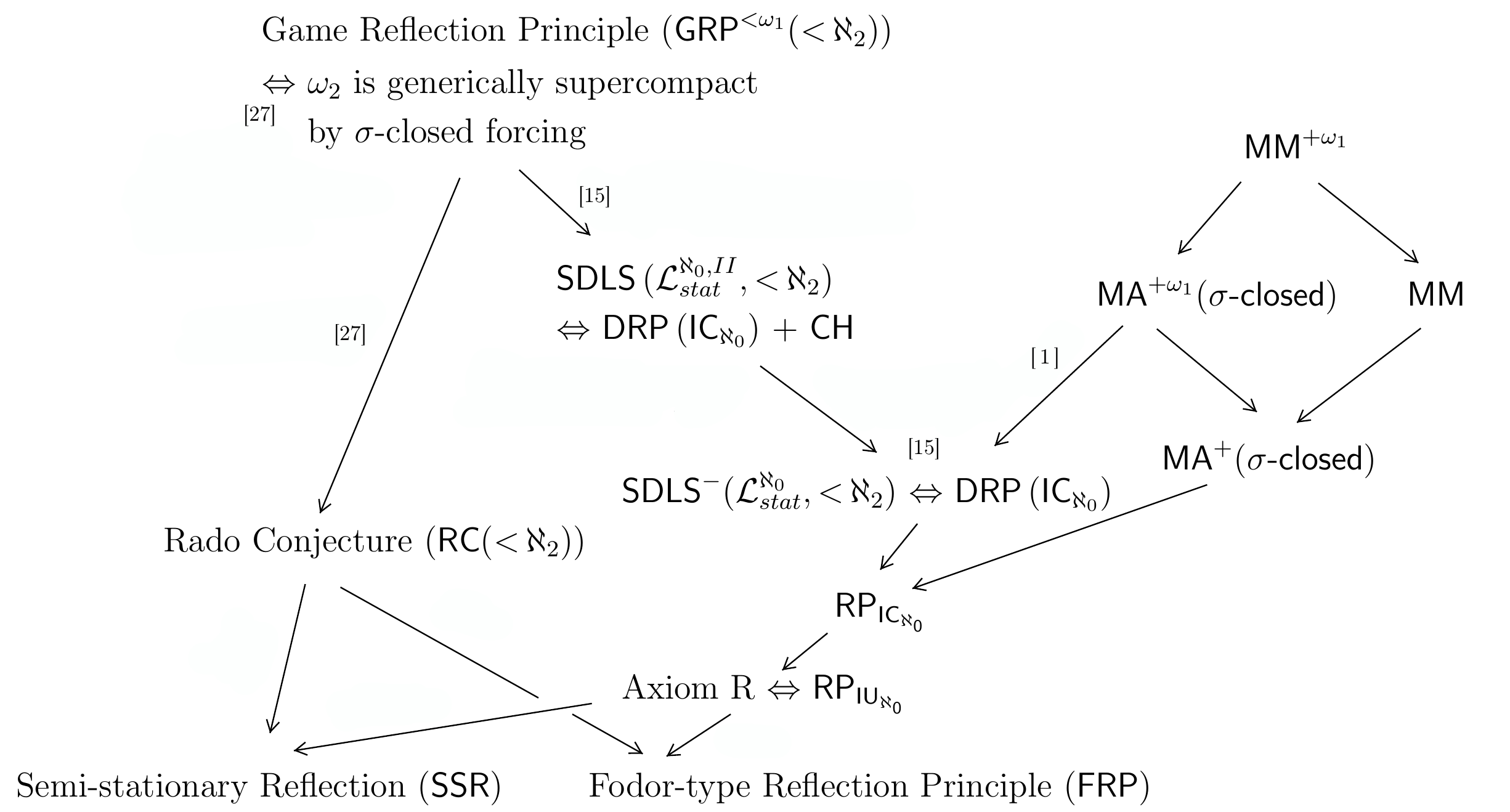}
\\[-1.8ex] 
\mbox{}\hspace{19em}{\footnotesize\rotatebox[origin=c]{45}{$\darkred\Leftrightarrow$}}\\[-2.4\jot]
\mbox{}\hspace{7em}{\darkred\footnotesize many ``mathematical'' reflection theorems with reflection down to
  $\darkred\LT\aleph_2$ }\\[-2\jot]
\mbox{}\hspace{7em}{\darkred\tiny{[4]},
  [5], [8], [14], [22], etc.}
\medskip\\
\phantomsection
\Label{figure1}
\centerline{{\tt Figure 1.}}
\bigskip

In the following we will elaborate on \xitemof{x-gen-4} above. 

\begin{Prop}
  \Label{p-gen-1} Suppose that $\kappa$ is $\calP$-gen.\ supercompact and $P$ is a property 
  of topological spaces which is \ifarxived\daimaru{a} \else ⓐ  \fi preserved by 
  homeomorphism and \ifarxived\daimaru{b} \else ⓑ  \fi downward absolute, meaning  
  that if $\uniW$ is a universe of set theory and $\uniW_0$ is an inner model in $\uniW$, 
  if a topological space $X\in\uniW$ satisfies the property $P$ in $\uniW$ then $X$ also 
  satisfies $P$ in $\uniW$. 
  Then, for any topological space $X$ of character $<\kappa$, if
  $\forces{\poP}{X\mbox{ satisfies }P}$ for any $\poP\in\calP$ then there is a 
  subspace $Y$ of $X$ of cardinality $\LT\kappa$ which satisfies $P$.\ifextended\else\qed\fi
\end{Prop}
{\ifextended\extendedcolor
\prf Suppose that $X$ is a topological space of character $\mu<\kappa$ and 
\begin{xitemize}
\ifarxived
\item[\daimaru{c}]
\else  
\item[ⓒ ]\fi $\forces{\poP}{X\mbox{ satisfies }P}$ for any $\poP\in\calP$.    
\end{xitemize}

\Wolog, we may assume that the underlying set of $X$ is a cardinal $\lambda\geq\kappa$ and 
the topology of $X$ is given by the system $\seqof{\tau_\alpha}{\alpha<\lambda}$ where  
$\tau_\alpha$ is an open nbhd basis of $\alpha$ of cardinality $\leq\mu$.

Let $\poP\in\calP$ be \st, for $(\uniV,\poP)$-generic $\genG$, there are $j$,
$M\subseteq\uniV[\genG]$ \st\quad\ifarxived\daimaru{1} \else ①  \fi
$\Elembed{j}{\uniV}{M}{\kappa}$,
\quad\ifarxived\daimaru{2} \else ②  \fi$j(\kappa)>\lambda$ and
\quad\ifarxived\daimaru{3} \else ③  \fi$j\imageof{\lambda}\in M$.

Let
\begin{xitemize}
\item[] $X':
  \begin{array}[t]{@{={}}l}
    \pairof{j\imageof{X},\,
      \setof{\pairof{\setof{j\imageof{U}}{U\in\tau_\alpha},j(\alpha)}}{\alpha\in\lambda}}\\
    \pairof{j\imageof{X},\,
      \setof{\pairof{\setof{V\cap j\imageof{X}}{V\in j(\tau_\alpha)},j(\alpha)}}{\alpha\in\lambda}}.
  \end{array}$
\end{xitemize}
Then $X'\in M$ by \ifarxived\daimaru{3} \else ③  \fi(and \ifarxived\daimaru{1}\else ①  \fi), 
$X'$ is a subspace of $j(X)$ and $X\cong X'$ (in
$\uniV[\genG]$). By \ifarxived\daimaru{a} \else ⓐ  \fi and\ 
\ifarxived\daimaru{c}\else ⓒ\,\fi, $\uniV[\genG]\modelof{X'\mbox{ satisfies }P}$ and hence
$M\modelof{X'\mbox{ satisfies }P}$ by \ifarxived\daimaru{b}\else ⓑ\,\fi. 
Thus, by \ifarxived\daimaru{2}\else ②\,\fi,  
$M\modelof{\mbox{there is a subspace of }X\mbox{ of size }\LT j(\kappa)\mbox{ satisfying }P}$
By elementarity, it follows that
  $V\modelof{\mbox{there is subspace of }X\mbox{ of size }\LT\kappa\mbox{ satisfying }P}$.
  \qedofProp\qedskip
\fi} 

\begin{Cor}\Label{p-gen-2} Suppose that \GRP\ holds. Then, for any topological 
  space $X$ of  
  character $\LE\aleph_1$ \st\ $\forces{\poP}{X\mbox{ is not metrizable}}$ for 
  any $\sigma$-closed \po\ $\poP$, there is a non-metrizable subspace of $X$ of cardinality
  $\LE\aleph_1$. 
\end{Cor}
\prf By \Thmof{p-gen-0} and \Propof{p-gen-1} for the property $P$ being non-metrizable. 
\qedofCor\qedskip

\section{Laver-generic large cardinals}
\Label{laver}
The Laver-genericity axioms (i.e.\ the axioms claiming the existence of Laver-generic large 
cardinals defined below) for respective classes of \pos complete the picture of reflection 
and absoluteness in terms of 
double plus versions of forcing axioms{\ifextended\extendedcolor\ given in \hyperref[figure1]{Figure 1} (see 
  \Thmof{p-laver-2} and 
Figure 3 below)\fi}.

A (definable) class $\calP$ of \pos\ is said to be {\It iterable} if\quad 
\ifextended
    {\extendedcolor \assertof{a} $\ssetof{\bbone}\in\calP$, \quad\assertof{b}
        $\calP$ is closed \wrt\ forcing 
        equivalence\quad (i.e.\ if $\poP\in\calP$ and $\poP\sim\poP'$ then $\poP'\in\calP$), 
        \quad\assertof{c} closed \wrt\ restriction\quad
        (i.e.\ if $\poP\in\calP$ then $\poP\restr\condp\in\calP$ for any $\condp\in\poP$),\quad and
        \quad\assertof{d} for 
        any $\poP\in\calP$ and $\poP$-name $\utpoQ$, $\forces{\poP}{\utpoQ\in\calP}$ implies
        $\poP\ast\utpoQ\in\calP$.}
\else
\assertof{a} $\calP$ is closed \wrt\ forcing 
equivalence\quad (i.e.\ if $\poP\in\calP$ and $\poP\sim\poP'$ then $\poP'\in\calP$), 
\quad\assertof{b} closed \wrt\ restriction\quad
(i.e.\ if $\poP\in\calP$ then $\poP\restr\condp\in\calP$ for any $\condp\in\poP$),\quad and\\
\assertof{c} for 
any $\poP\in\calP$ and $\poP$-name $\utpoQ$, $\forces{\poP}{\utpoQ\in\calP}$ implies
$\poP\ast\utpoQ\in\calP$.\fi

For an iterable class $\calP$ of \pos, 
a cardinal $\kappa$ is said to be {\It $\darkred\calP$-Laver-gen.\ supercompact} 
\footnote{The definition of  
  Laver-generic large cardinals given here is slightly stronger than the one given in 
  \cite{sfetal-II}.  The Laver-generic large cardinals in the sense of present subsection is 
  called {\it strongly} Laver-generic large cardinals in \cite{sfetal-II}.} 
if, for 
any $\lambda\geq\kappa$ and $\poP\in\calP$, there 
is a $\poP$-name $\utpoQ$ with
$\forces{\poP}{\utpoQ\in\calP}$ \st, for 
$(\uniV,\poP\ast\utpoQ)$-generic 
$\genH$, there are $j$, $M\subseteq\uniV[\genH ]$ with
\quad\ifarxived\daimaru{1} \else ①  \fi$\Elembed{j}{\uniV}{M}{\kappa}$,\\
\ifarxived\daimaru{2} \else②  \fi$j(\kappa)>\lambda$, and
\quad\ifarxived\daimaru{3}' \else ③' \fi$\poP\ast\utpoQ$, $\genH$,
$j\imageof{\lambda}\in M$. 

Recall that a cardinal $\kappa$ is {\it\darkred superhuge} ({\it\darkblue 
  super-almost-huge}, resp.) if, for any $\lambda>\kappa$, 
there are classes $j$, $M$ \st\quad\ifarxived\daimaru{1} \else ①\ \fi
$\Elembed{j}{\uniV}{M}{\kappa}$,\quad\ifarxived\daimaru{2} \else ② \fi
$j(\kappa)>\lambda$ and\quad
\ifarxived\\\daimaru{3} \else ③  \fi$\darkred\fnsp{j(\kappa)}{M}\subseteq M$\quad
({$\darkblue\fnsp{j(\kappa)\GT}{M}\subseteq M$}, resp.).

These notions of large cardinals can be straightforwardly translated into their 
Laver-generic versions: For an iterable class $\calP$ of \pos, $\kappa$ is 
{\it\darkred $\darkred\calP$-Laver-gen.\ superhuge} ({\it\darkblue $\darkblue\calP$-Laver-gen.\ 
  super-almost-huge}, resp.) if, for 
any $\lambda\geq\kappa$, $\poP\in\calP$, there 
is a $\poP$-name $\utpoQ$ with
$\forces{\poP}{\utpoQ\in\calP}$ \st, for 
$(\uniV,\poP\ast\utpoQ)$-generic 
$\genH$, there are $j$, $M\subseteq\uniV[\genH ]$ with\quad 
\ifarxived\daimaru{1} \else ① \fi$\Elembed{j}{\uniV}{M}{\kappa}$,\quad 
\ifarxived\daimaru{2} \else ② \fi$j(\kappa)>\lambda$, and\quad 
\ifarxived\daimaru{3}' \else ③' \fi$\poP$, $\poP\ast\utpoQ$,
$\genH\in M$, and {$\darkred j\imageof{j(\kappa)}\in M 
  $} ({\darkblue $j\imageof\mu\in M$ for all $\darkblue\mu<j(\kappa)$}, resp.).

Sometimes it is more convenient to consider the following additional property
which we called the tightness of Laver-genericity: For an iterable $\calP$, 
a $\calP$-Laver-gen.\ supercompact cardinal ($\calP$-Laver-gen.\ huge cardinal, etc., 
resp.) is {\It tightly $\darkred\calP$-Laver-gen.\  
  supercompact} ({\it\darkblue tightly $\calP$-Laver-gen.\ huge},  etc., resp.) if 
  the condition
  \begin{xitemize}
  \ifarxived
  \item[\daimaru{4}]
  \else
  \item[④ ] \fi
    $\poP\ast\utpoQ$ is forcing equivalent to a \po\ of 
    cardinality $\leq j(\kappa)$.
  \end{xitemize}
additionally holds for the elementary embedding $j$ in the definition. 

The strongest notion of large cardinal we consider in this paper in connection with its 
Laver-generic version is that of ultrahuge cardinal introduced by
Tsaprounis 
\cite{tsaprounis2}. A cardinal $\kappa$ is {\It ultrahuge} if for any $\lambda>\kappa$ there is
  $\Elembed{j}{\uniV}{M}{\kappa}$ \st\ $j(\kappa)>\lambda$ and
$\fnsp{j(\kappa)}{M}, V_{j(\lambda)}\subseteq M$. In terms of consistency strength 
ultrahuge cardinal is placed between superhuge and 2-almost-huge (Theorem 3.4 in 
\cite{tsaprounis2}).

For an iterable class $\calP$ of \pos, a cardinal $\kappa$ is ({\it\darkblue 
  tightly}) {\it\darkred $\darkred \calP$-Laver gen.\ ultrahuge}, if, for 
  any $\lambda>\kappa$ and $\poP\in\calP$ there is a $\poP$-name $\utpoQ$ with 
$\forces{\poP}{\utpoQ\in\calP}$, \st\  
  for $(\uniV,\poP\ast\utpoQ)$-generic $\genH$, there  
  are $j, M\subseteq\uniV[\genH]$ 
  \st\ $\Elembed{j}{\uniV}{M}{\kappa}$, $j(\kappa)>\lambda$, 
  $\poP,\genH,(V_{j(\lambda)})^{\uniV[\genH]}\in M$
  ({\darkblue and
    $\darkblue\darkblue\poP\ast\utpoQ$ is forcing 
    equivalent to a \po\ of size $\darkblue j(\kappa)$}).

The following theorem is used to construct models with a Laver-generically ultrahuge 
cardinal: 
{\ifextended\extendedcolor
\begin{ThmA}\Label{p-gen-3}\extendedcolor{\rm (Tsaprounis \cite{tsaprounis2})} If
  $\kappa$ is an ultrahuge cardinal, then $\kappa$ carries a Laver function
  $\mapping{f}{\kappa}{V_\kappa}$, i.e. a function $f$ with the property:
  \begin{xitemize}
  \xitem[x-gen-5] for every cardinal $\lambda\geq\kappa$ and any $x\in\calH(\lambda^+)$ 
    there is an $\Elembed{j}{\uniV}{M}{\kappa}$ with $j(\kappa)>\lambda$, and
    $\fnsp{j(\kappa)}{M}, V_{j(\lambda)}\subseteq M$  \st\ $x=j(f)(\kappa)$. 
    \qed\vspace{-1.44ex}
  \end{xitemize}
\end{ThmA}\fi} 

By definition,  it is obvious that we have the following implications: \bigskip\bigskip

\noindent\hspace{-0.5em}
\scalebox{0.53}[0.65]{
$\begin{array}[t]{ccccccccc}
  \vbox{\hbox{tightly $\calP$-Laver-gen.} 
    \hbox{ultrahuge}}
  &\Rightarrow
  &\vbox{\hbox{tightly $\calP$-Laver-gen.} 
    \hbox{superhuge}}
  &\Rightarrow
  &\vbox{\hbox{tightly $\calP$-Laver-gen.}
    \hbox{super-almost-huge}}
  &\Rightarrow
  &\vbox{\hbox{tightly $\calP$-Laver-gen.}
    \hbox{supercompact}}
  &\Rightarrow 
  &\vbox{\hbox{tightly $\calP$-Laver-gen.}
    \hbox{measurable}}
  \\[-0.8\jot]
  \scalebox{1.2}[0.9]{$\Downarrow$} && 
  \scalebox{1.2}[0.9]{$\Downarrow$} && \scalebox{1.2}[0.9]{$\Downarrow$}
  && \scalebox{1.2}[0.9]{$\Downarrow$}  && \scalebox{1.2}[0.9]{$\Downarrow$}\\[-1\jot]
  \calP\mbox{-Laver-gen.\ ultrahuge} &\Rightarrow 
  &\calP\mbox{-Laver-gen.\ superhuge} &\Rightarrow &\calP\mbox{-Laver-gen.\ super-almost-huge}
  &\Rightarrow &\calP\mbox{-Laver-gen.\ supercompact}
  &\Rightarrow &\calP\mbox{-Laver-gen.\ measurable}
  \\[-0.8\jot]
  \scalebox{1.2}[0.9]{$\Downarrow$} &&
  \scalebox{1.2}[0.9]{$\Downarrow$} && \scalebox{1.2}[0.9]{$\Downarrow$}
  && \scalebox{1.2}[0.9]{$\Downarrow$}  && \scalebox{1.2}[0.9]{$\Downarrow$}\\[-1\jot]
  \calP\mbox{-gen.\ ultrahuge} &\Rightarrow 
  &\calP\mbox{-gen.\ superhuge} &\Rightarrow &\calP\mbox{-gen.\ super-almost-huge}
  &\Rightarrow &\calP\mbox{-gen.\ supercompact}
  &\Rightarrow &\calP\mbox{-gen.\ measurable}
\end{array}
  $}\bigskip\\
\phantomsection
\Label{figure2}
\centerline{{\tt Figure 2.}}
\bigskip

Some of the horizontal implications should be irreversible. At the moment however we can only prove 
the irreversibility of the implication from (tightly) $\calP$-(Laver)-gen.\ ultrahugeness 
to (tightly) $\calP$-(Laver)-gen.\ supercompactness. 

\begin{Prop}{\ifextended\extendedcolor\rm (Proposition 4 in \cite{slides-x})\fi}
  \Label{p-laver-a} Suppose that $\poP$ is a class of \pos\ \st\ there is a construction of 
  a model with a tightly $\calP$-Laver gen.\ supercompact cardinal starting from an 
  arbitrary model with an supercompact cardinal $\kappa$ by a \po\ of 
  cardinality $\kappa$.\footnotemark 
  Then tightly $\calP$-Laver gen.\ supercompactness of $\kappa$ 
  does not necessarily imply the $\calP$-gen.\ super-almost-hugeness. 
\end{Prop}
\footnotetext{By the following \Thmof{p-laver-0}, the class of all $\sigma$-closed \pos\ 
  and the class of all ccc \pos\ satisfy this condition.}

For the proof of \Propof{p-laver-a} we use the following observation:
\begin{Lemma}\Label{p-laver-a-0}
  Suppose that $\kappa$ is $\calP$-gen. ultrahuge   
  for an arbitrary class $\calP$ of \pos. 
  If there is an inaccessible $\lambda_0>\kappa$ then there are cofinally many inaccessible 
  in $\uniV$.
\end{Lemma}
\prf 
Let $\lambda>\lambda_0$ be an arbitrary cardinal. Then 
there is $\poP\in\calP$ \st, for $(\uniV,\poP)$-generic $\genG$, there are 
$j, M\subseteq\uniV[\genG]$ \st\ 
\begin{xitemize}
  \xitem[x-laver-1] 
    $\Elembed{j}{\uniV}{M}{\kappa}$,\qquad
  \ixitem[x-laver-2] 
    $j(\kappa)>\lambda$, 
    and
  \xitem[x-laver-3] 
    $(V_{j(\lambda)})^{\uniV[\genG]}\in M$. 
\end{xitemize}

By \xitemof{x-laver-2} and elementarity
\xitemof{x-laver-1}, 
we have $j(\lambda_0)>\lambda$.
By elementarity 
\xitemof{x-laver-1}, $M\modelof{j(\lambda_0)\mbox{ is inaccessible}}$. 
By
\xitemof{x-laver-3}, $\uniV[\genG]\modelof{j(\lambda_0)\mbox{ is inaccessible}}$, and hence
$\uniV\modelof{j(\lambda_0)\mbox{ is inaccessible}}$.
\qedofLemma\qedskip

\noindent{\bf Proof of \bfPropof{p-laver-a}:} 
Suppose that $\kappa$ is a supercompact cardinal 
and $\lambda_0>\kappa$ is an inaccessible cardinal. 

We may assume that $\lambda_0$ is the 
largest inaccessible cardinal: if there is inaccessible cardinal larger than $\lambda_0$, 
then let $\lambda_1$ be the least such inaccessible cardinal. In $V_{\lambda_1}$, 
$\lambda_0$ is the largest inaccessible cardinal and $\kappa$ is supercompact (see e.g.\ 
Exercise 22.8,\,\assertof{a} in \cite{higher-inf}
{\ifextended\extendedcolor 
  : $V_{\lambda_1}\modelof{\kappa\mbox{ is supercompact}}$ can be seen using the 
  characterization of supercompactness in terms of ultrafilters.\fi}).

Let $\poP$ be a \po\ of size $\kappa$ \st, for
$(\uniV,\poP)$-generic $\genG$, we have
$\uniV[\genG]\modelof{\kappa\xmbox{ is tightly }\poP\xmbox{-Laver gen.\ supercompact}}$.
Note that 
$\uniV[\genG]\modelof{\lambda_0\xmbox{ is the largest inaccessible cardinal}}$. 
Thus, by \Lemmaof{p-laver-a-0}, it follows that
$\uniV[\genG]\modelof{\kappa\xmbox{ is not }\calP\xmbox{-gen.\ ultrahuge}}$.
\qedof{\Propof{p-laver-a}} 
\qedskip

(Tightly) Laver-generic large cardinal is actually first-order definable (i.e. it has a 
characterization formalizable in the language of \ZFC\,), cf.\ \cite{fuchino-sakai}.
Thus ``Forcing Theorems'' are available for arguments with Laver-genericity. 
Because of this and because an iterable class $\calP$ is closed under 
restriction to a condition, by definition, we may be lazy about the quantification on 
generic filters like in the context of ``\,for 
a/any $(\uniV,\poP\ast\utpoQ)$-generic $\genH$ ...''

The Examples \ref{ex-gen-1}, \ref{ex-gen-2} are actually examples of the construction of 
models with a Laver-generic large cardinal.
\begin{Thm}{\rm (Theorem 5.2, \cite{sfetal-II})} 
  \Label{p-laver-0}  \wassertof{1} Suppose that $\kappa$ is {\darkred supercompact}
  ({\darkblue superhuge}, etc., resp.) and $\poP=\Col(\aleph_1,\kappa)$. Then, 
  in $\uniV[\genG]$, for any $(\uniV,\poP)$-generic $\genG$, $\aleph_2^{\uniV[\genG]}$
  ($=\kappa$) is tightly 
  $\sigma$-closed-Laver-gen.\ {\darkred supercompact} (
    {\darkblue superhuge}, etc., resp.) 
    and $\CH$ holds.\smallskip

  \wassert{2}
  Suppose that $\kappa$ is {\darkgreen super-almost-huge}  ({\darkblue superhuge}, etc., resp.) with a
  Laver function $\mapping{f}{\kappa}{V_\kappa}$ for {\darkgreen super-almost-hugeness} 
  ({\darkblue superhugeness},  
  etc., resp.), and $\poP$ is the CS-iteration for forcing \PFA\ along with $f$. Then,  
  in $\uniV[\genG]$ for any $(\uniV,\poP)$-generic $\genG$, $\aleph_2^{\uniV[\genG]}$
  ($=\kappa$)  
  is tightly proper-Laver-gen.\ {\darkgreen super-almost-huge} ({\darkblue superhuge}, etc.,
  resp.) and 
  $\continuum=\aleph_2$ holds.\footnotemark\smallskip

  \wassert{2\rlap{$'$}\,}
  Suppose that $\kappa$ is {\darkgreen super-almost-huge}  ({\darkblue superhuge}, resp.) with a
  Laver function $\mapping{f}{\kappa}{V_\kappa}$ for {\darkgreen super-almost-hugeness} 
  ({\darkblue superhugeness},  etc.,
  resp.), and $\poP$ is the RCS-iteration for forcing \MM\ along with $f$. Then,  
  in $\uniV[\genG]$ for any $(\uniV,\poP)$-generic $\genG$, $\aleph_2^{\uniV[\genG]}$
  ($=\kappa$)  
  is tightly semi-proper-Laver-gen.\ {\darkgreen super-almost-huge} ({\darkblue superhuge}, etc.,
  resp.) and 
  $\continuum=\aleph_2$ holds.\addtocounter{footnote}{-1}\footnotemark\smallskip

  \wassert{3}
  Suppose that $\kappa$ is {\darkred supercompact} (
  {\darkblue superhuge}, etc., resp.) with a 
  Laver function $\mapping{f}{\kappa}{V_\kappa}$ for {\darkred supercompactness} (
  {\darkblue superhugeness}, etc., resp.), and $\poP$ is a FS-iteration for forcing \MA\ 
  along with $f$.  
  Then, in $\uniV[\genG]$ for any $(\uniV,\poP)$-generic $\genG$, $\continuum$ ($=\kappa$) is 
  tightly ccc-Laver-gen.\ {\darkred supercompact} (
  {\darkblue superhuge}, etc., resp.).  $\kappa=\continuum$, and $\kappa$ is 
  very large.\ifextended\else\qed\fi
\end{Thm}
\footnotetext{It seems that the construction does not work with supercompact $\kappa$ here.}
{\ifextended\extendedcolor

In the following we give a proof of the case \assertof{2} of \Thmof{p-laver-0}
for ultrahugeness and its 
Laver-generic version. For this case,  we need the next lemma. Note that $\kappa$ is almost 
huge if it is ultrahuge.
\begin{Lemma}\Label{p-laver-0-0}
  \wassertof{1} If $\kappa$ is almost huge and $\Elembed{j}{\uniV}{M}{\kappa}$ is 
  an almost huge elementary embedding, then the target of $j$ (i.e. $j(\kappa)$) is 
  inaccessible.\smallskip 

  \wassert{2} Suppose that $\kappa$ is super almost huge. Then there are cofinally many 
  inaccessible cardinals.  
\end{Lemma}
\prf \assertof{1}: Suppose that $\Elembed{j}{\uniV}{M}{\kappa}$ is 
an almost huge elementary embedding. Thus we have in particular \ixitem[x-laver-4]
$\fnsp{j(\kappa)\GT}{M}\subseteq M$. Since $\kappa$ is inaccessible,
$M\modelof{j(\kappa)\mbox{ is inaccessible}}$ by elementarity. By \xitemof{x-laver-4}, it 
follows that $j(\kappa)$ is really inaccessible. \smallskip

\assertof{2}: follows from \assertof{1} since, if $\kappa$ is super almost huge, the 
targets of almost huge elementary embeddings are unbounded. 
\memo{kobe2023-05-29-pf\\Lemma 12}
\qedofLemma\qedskip

\noindent
\prfof{\bfThmof{p-laver-0}} 
We prove \assertof{2} 
for ultrahugeness and its 
Laver-generic version. Other cases can be proved similarly. 

Suppose that $\kappa$ is ultrahuge and $\mapping{f}{\kappa}{V_\kappa}$ is an 
ultrahuge Laver function (see \ThmAof{p-gen-3}).

Let $\seqof{\poP_\alpha,\utpoQ_\beta}{\alpha\leq\kappa,\beta<\kappa}$ be a 
CS-iteration with 

\begin{xitemize}
\item[] 
  $\utpoQ_\beta:=\left\{\,
  \begin{array}{@{}ll}
    f(\beta), &\mbox{if }\forces{\poP_\beta}{f(\beta)\mbox{ is a proper \po};}\\
    \bbone, &\mbox{otherwise.}
  \end{array}
  \right.$
\end{xitemize}

We show that $\poP_\kappa$ forces that $\kappa$ is tightly proper-Laver generically 
ultrahuge. 

Let $\genG_\kappa$ be a $(\uniV,\poP_\kappa)$-generic filter. 
Suppose $\lambda>\kappa$ and $\poP$ be a proper \po\  
in $\uniV[\genG_\kappa]$. Let $\utpoP$ be a $\poP_\kappa$-name of $\poP$. 

By \Lemmaof{p-laver-0-0},\,\assertof{2}, we may assume that $\lambda$ is inaccessible.
\smallskip

Let $\Elembed{j}{\uniV}{M}{\kappa}$ be \st\ $j(f)(\kappa)=\utpoP$,
$j(\kappa)>\lambda$, and\quad
\begin{xitemize}
\xitem[x-laver-5]
  $\fnsp{j(\kappa)}{M}, V_{j(\lambda+1)}\subseteq M$.\footnotemark
\end{xitemize}
\footnotetext{It is intentional that we 
  choose $M$ here with slightly stronger closure property by saying
  $V_{j(\lambda+1)}\subseteq M$ instead of $V_{j(\lambda)}\subseteq M$.}

By elementarity, we have

\begin{xitemize}
\item[] 
  $M\modelof{\,
    \begin{array}[t]{@{}l}
      j(\poP_\kappa)\mbox{ is a CS-iteration }
      \seqof{\poP^*_\alpha,\utpoQ^*_\beta}{\alpha\leq j(\kappa),\beta<j(\kappa)}
      \mbox{ of proper }\\
      \mbox{\pos\ with the book-keeping }j(f)\mbox{ and }
      \cardof{\poP^*_\alpha}<j(\kappa)\mbox{ for 
        all }\alpha<\kappa}.
    \end{array}$
\end{xitemize}

Note that $\poP^*_\alpha=\poP_\alpha$ for all $\alpha\leq\kappa$, $\poP_\kappa\in M$, and
$\utpoQ^*_\kappa=\utpoP$.
Thus, by the Factor Lemma

\begin{xitemize}
\item[] 
  $M[\genG_\kappa]\modelof{\,
  \begin{array}[t]{@{}l}
    j(\poP_\kappa)/\genG_\kappa\mbox{ is (forcing equivalent to) a CS-iteration of proper}\\
    \mbox{\pos\ of length }j(\kappa)\mbox{ and its 0th iterand is }\poP\,}.
  \end{array}$
\end{xitemize}
By the $\kappa$-cc of $\poP_\kappa$ and by \xitemof{x-laver-5},
we have
  $\fnsp{\lambda}{(M[\genG_\kappa])}\subseteq M[\genG_\kappa]$. 
\begin{xitemize}
\item[] 
  $\uniV[\genG_\kappa]\modelof{\,
    \begin{array}[t]{@{}l}
      j(\poP_\kappa)/\genG_\kappa\rlap{\mbox{ is (forcing equivalent to) a CS-iteration of proper}}\\
      \mbox{\pos\ of length }j(\kappa)\mbox{ and its 0th iterand is }\poP\,}.
    \end{array}$
\end{xitemize}

It follows that, in $\uniV[\genG_\kappa]$, we have
$j(\poP_\kappa)/\genG_\kappa\sim \poP\ast\utpoQ^*$  
where
\begin{xitemize}
\item[] 
  $\uniV[\genG_\kappa]\models{\forces{\poP}{\utpoQ^*\mbox{ is proper}}}$.
\end{xitemize}

Let $\genH$ be a $(\uniV[\genG_\kappa], j(\poP_\kappa)/\genG_\kappa)$-generic 
filter: Note that $\genH$ corresponds to a
$(\uniV[\genG], \poP\ast\utpoQ^*)$-generic filter, and $\genG_\kappa\ast\genH$ 
corresponds to a $(\uniV,j(\poP_\kappa))$-generic filter extending $\genG_\kappa$. I shall 
denote the latter also with $\genG\ast\genH$.

Let $\tilde{j}$ be the ``lifting'' of $j$ defined by
\begin{xitemize}
\item[] 
  $\mapping{\tilde{j}}{\uniV[\genG_\kappa]}{M[\genG_\kappa\ast\genH]};$\quad
  $\utilde{a}[\genG_\kappa]\ \mapsto\ j(\utilde{a})[\genG_\kappa\ast\genH]$\quad 
  for all $\poP_\kappa$-name $\utilde{a}$.
\end{xitemize}

Then we have $j\subseteq\tilde{j}$,\quad
$\Elembed{\tilde{j}}{\uniV[\genG_\kappa]}{M[\genG_\kappa\ast\genH]}{\kappa}$,
$\tilde{j}\imageof{\lambda}= j\imageof{\lambda}\in M\subseteq M[\genG_\kappa\ast\genH]$,
$\cardof{j(\poP_\kappa)/\genG_\kappa}^{\uniV[\genG_\kappa]}
\leq\cardof{j(\poP_\kappa)}^M=j(\kappa)$.

$\genG_\kappa\ast\genH$ seen as a $(\uniV,j(\poP_\kappa))$-gen.\ filter has 
cardinality $j(\kappa)<j(\lambda)$ and it is $\in V_{j(\lambda)}$.

Thus, there is a $j(\poP_\kappa\ast\utpoQ)$-name $\utilde{V}$ of 
$(V_{j(\lambda)})^{\uniV[\genG_\kappa\ast\genH]}$ in $V_{j(\lambda)+1}=V_{j(\lambda+1)}$.
  
It follows that 
\begin{xitemize}
\item[] 
  $(V_{j(\lambda)})^{\uniV[\genG_\kappa\ast\genH]}
  =\utilde{V}[\genG_\kappa\ast\genH]\in M[\genG_\kappa\ast\genH]$.
\end{xitemize}

This shows that
$\uniV[\genG_\kappa]\modelof{\kappa\mbox{ is tightly proper-Laver-gen.\ ultrahuge}}$. \\
\qedof{\Thmof{p-laver-0}}\qedskip

\fi}

The circumstance that the three possibilities of the cardinality of the 
continuum: $\aleph_1$, $\aleph_2$, or 
very large, are highlighted in \Thmof{p-laver-0}, has also an explanation in terms of 
Laver-genericity:  
\begin{Thm}{\rm (The Trichotomy Theorem \cite{sfetal-II}, see also \cite{nagoya})}
  \Label{p-laver-1} 
  \wassertof{A} If $\kappa$ is ${\calP}$-Laver-gen.\ supercompact for an iterable 
  class ${\calP}$ of  
  \pos\  such that \assertof{a} \ 
  all $\poP\in{\calP}$ are $\omega_1$ preserving, 
  \assertof{b} \  all $\poP\in{\calP}$ do not add reals, and \assertof{c} \  there is a
  $\poP_1\in{\calP}$ which 
  collapses $\omega_2$, \underline{then} $\kappa=\aleph_2$ and \CH\ holds.\smallskip

  \wassert{B} If $\kappa$ is ${\calP}$-Laver-gen.\ supercompact for an iterable class ${\calP}$ of 
  \pos\  such that \assertof{a} \   all $\poP\in{\calP}$ are $\omega_1$-preserving, 
  \assertof{b$'$}\ there is  
  a $\poP_0\in{\calP}$ which add a real, and \assertof{c} \  there is a $\poP_1$ which 
  collapses $\omega_2$, \underline{then} $\kappa=\aleph_2\leq 2^{\aleph_0}$. If $\calP$ 
  contains enough many proper \pos\ then $\kappa=\aleph_2=2^{\aleph_0}$ (For 
    the last assertion see \Thmof{p-laver-2} below). \smallskip

  \wassert{$\Gamma$}
  If $\kappa$ is ${\calP}$-Laver-gen.\ supercompact for an iterable class ${\calP}$ of 
  \pos\  such that \assertof{a$'$}\ \ all $\poP\in{\calP}$ preserve cardinals, and 
  \assertof{b$'$}\ \ there is a 
  $\poP_0\in{\calP}$ which  
  adds a real, \underline{then} $\kappa$ is ``very large'' and
  $\kappa\leq 2^{\aleph_0}$. 
  If $\kappa$ is tightly $\calP$-Laver-gen.\ superhuge then $\kappa=2^{\aleph_0}$. 
  \ifextended\else\qed\fi 
\end{Thm}
{\ifextended\extendedcolor 

\Thmof{p-laver-1} follows from the next \Lemmaof{p-Lg-RA-1-1}, \Lemmaof{p-Lg-RA-1-2-0}, and 
\Thmof{p-MP-2-0}.\qedskip

\begin{Lemma}\Label{p-Lg-RA-1-1}
  Suppose that $\calP$ is a class of \pos\ and $\kappa$ is $\calP$-generically measurable. 
  Then:\smallskip

  \wassert{1} $\kappa$ is regular.\smallskip

  \wassert{2} If all elements of $\calP$ preserve $\kappa$, then $\kappa$  
  is a weakly inaccessible cardinal.\smallskip

  \wassert{3} If all elements of $\calP$ also preserve regularity of $\kappa$, 
  then $\kappa$ is a 
  weakly inaccessible cardinal which is a stationary limit of weakly inaccessible, limit of 
  limits of weakly inaccessible cardinals, etc.  \smallskip

  \wassert{4} If all elements of $\calP$ preserve stationarity of subsets of all regular
  $\lambda\leq\kappa$, then 
  $\kappa$ is a weakly 
  inaccessible which is a stationary limit of weakly inaccessible, stationary limit of 
  stationary limits of weakly inaccessible cardinals, etc. 
\end{Lemma}
\prf Let $\poP\in\calP$ be \st, for $(\uniV,\poP)$-generic 
$\genG$, there are $j$, $M\subseteq\uniV[\genG]$ \st\ $\Elembed{j}{\uniV}{M}{\kappa}$.

\assertof{1}: If $\kappa$ were not regular, there would be $\mu<\kappa$ and 
cofinal $\mapping{f}{\mu}{\kappa}$. $j(f)=f$ by $\crit(j)=\kappa$. Hence, by elementarity,
$M\models j(\kappa)=\sup_{\alpha<\mu}f(\alpha)=\kappa$. This is a contradiction.\smallskip

\assertof{2}: Suppose that $\poP$ preserves $\kappa$.
$\kappa$ is a limit cardinal:
If $\kappa$ were a successor cardinal, there would be $\mu<\kappa$ with a sequence
$\seqof{f_\alpha}{\alpha<\kappa}$ \st\ each $f_\alpha$ is a surjection form $\mu$ to
$\alpha$ ($\mu$ to $\alpha+1$ for finite $\alpha$). Let
$\seqof{f^*_\alpha}{\alpha<j(\kappa)}:=j(\seqof{f_\alpha}{\alpha<\kappa})$. Then $f^*_\kappa$ 
is a surjection from $\mu$ to $\kappa$ (in $M$, and hence in $\uniV[\genG]$). This is a 
contradiction to $\poP\in\calP$.

Together with \assertof{1}, this implies that $\kappa$ is a weakly inaccessible cardinal.
\smallskip

\assertof{3}: Let $\poP\in\calP$ be \st\ for a $(\uniV,\poP)$-generic $\genG$ there 
are $j$, $N\subseteq\uniV[\genG]$ \st\ $\Elembed{j}{\uniV}{M}{\kappa}$. $\kappa$ is weakly 
inaccessible by \assertof{2}.

Suppose that $D\subseteq\kappa$ is a club (in $\uniV$). We want to show that $D$ contains a 
weakly inaccessible cardinal. Since $j(D)$ is closed by elementarity and since $j(D)\cap\kappa= D$ ,
$\kappa\in j(D)$. Since $\poP$ preserves confinality,
$\uniV[\genG]\modelof{\kappa\mbox{ is weakly inaccessible}}$. Hence
$M\modelof{\kappa\mbox{ is weakly inaccessible}}$. Thus
$M\modelof{j(D)\mbox{ contains a weakly inaccessible cardinal}}$. By elementarity, it 
follows that $\uniV\modelof{D\mbox{ contains a weakly inaccessible cardinal}}$. 
\smallskip

\assertof{4}: Let $\poP$, $\genG$, $j$, $M$ be as in \assertof{3} where $\poP$ now 
preserves stationary subsets of $\kappa$.

Suppose that $D\subseteq\kappa$ is a club (in $\uniV$). We want to show that $D$ contains a 
weakly inaccessible cardinal which is a stationary limit of weakly inaccessible cardinals.

By \assertof{3}, $S=\setof{\alpha<\kappa}{\alpha\mbox{ is weakly inaccessible}}$ is 
stationary. It follows that $S$ remains stationary subset of $\kappa$ in $\uniV[\genG]$. 
Since $S=j(S)\cap\kappa$, We have $S\in M$. As in \assertof{3}, we also have
$\kappa\in j(D)$. 

Thus
$M\modelof{\mbox{there is a weakly inaccessible }\delta\in j(D)\xmbox{ which is a stationary limit 
of inaccessible cardinals}}$. By elementarity it follows that 
$\uniV\modelof{\mbox{there is a weakly inaccessible}\delta\in D\mbox{ which is a stationary limit 
of inaccessible cardinals}}$. 
\qedofLemma
\qedskip

\begin{Lemma}\Label{p-Lg-RA-1-2-0}{\rm(Proposition 4, in \cite{nagoya})} \wassertof{1} If
  $\kappa$ is generically measurable by $\omega_1$ preserving 
  $\calP$, then $\omega_1<\kappa$.\smallskip

  \wassert{2} If $\kappa$ is Laver-gen.\ supercompact for $\omega_1$-preserving $\calP$ with
  $\Col(\omega_1,\ssetof{\omega_2})\in\calP$ then $\kappa=\omega_2$.\smallskip

  \wassert{3} If $\kappa$ is Laver-gen.\ supercompact for $\calP$ which contains a \po\ 
  adding a new real then $\kappa\leq2^{\aleph_0}$.\smallskip

  \wassert{4} If $\kappa$ is gen.\ supercompact by $\calP$ \st\ all \pos\ in $\calP$ do not 
  add any reals then $2^{\aleph_0}<\kappa$.
\end{Lemma}
\prf \assertof{1}: Suppose $\kappa\leq\omega_1$. Since $\kappa=\omega$ is impossible, we 
have $\kappa=\omega_1$. Let $\poP\in\calP$ and $\genG$ 
be $(\uniV,\poP)$-generic \st\ $\Elembed{j}{\uniV}{M\subseteq\uniV[\genG]}{\kappa}$.
Then we have $M\modelof{j(\kappa)=\omega_1}$ by elementarity. Since $\kappa<j(\kappa)$, 
$M\modelof{\kappa\mbox{ is countable}}$. Thus
$\uniV[\genG]\modelof{\kappa\mbox{ is countable}}$. This is a contradiction to the 
assumption that elements of $\calP$ are $\omega_1$ preserving. 
\smallskip

\assertof{2}: Suppose that $\kappa\not=\omega_2$. By \assertof{1}, we then have
\begin{xitemize}
\xitem[x-Lg-0] 
  $\kappa>\omega_2$. 
\end{xitemize}

Let $\poP:=\Col(\omega_1,\ssetof{\omega_2})$, and let $\utpoQ$ be 
a $\poP$-name with $\forces{\poP}{\utpoQ\in\calP}$ \st\ there is 
a $(V,\poP\ast\utpoQ)$-generic $\genH$ with $j$, $M\subseteq\uniV[\genH]$ \st\
$\Elembed{j}{\uniV}{M}{\kappa}$, and $\genH\in M$. By \xitemof{x-Lg-0},
$j({\omega_2}^\uniV)={\omega_2}^\uniV$ and $M\modelof{j({\omega_2}^\uniV)=\omega_2}$ by 
elementarity. On the other hand, $\genH$ codes a collapsing of ${\omega_2}^\uniV$. Thus
$M\modelof{\cardof{{\omega_2}^\uniV}=\aleph_1}$. This is a contradiction. \smallskip

\assertof{3}: Suppose that $\mu<\kappa$ and $\seqof{a_\alpha}{\alpha<\mu}$ is a sequence of 
reals. We show that $\seqof{a_\alpha}{\alpha<\mu}$ is not an enumeration 
of $\psof{\omega}$.

Let $\poP\in\calP$ be a \po\ which adds a real, and let $\utpoQ$ be a $\poP$-name of a \po\ 
\st, for a $(\uniV,\poP\ast\utpoQ)$-generic $\genH$ there are $j$, $M\subseteq\uniV[\genH]$ 
\st\ $\Elembed{j}{\uniV}{M}{\kappa}$ and $\genH\in M$. By elementarity and 
since $\mu<\kappa$, we have $j(\seqof{a_\alpha}{\alpha<\mu})=\seqof{a_\alpha}{\alpha<\mu}$. 
Since there is a new reals coded in $\genH$, 
\begin{xitemize}
\item[] $M\modelof{\seqof{a_\alpha}{\alpha<\mu}\mbox{ is not an enumeration of }\psof{\omega}}$.
\end{xitemize}
By elementarity it follows that
\begin{xitemize}
\item[] $\uniV\modelof{\seqof{a_\alpha}{\alpha<\mu}\mbox{ is not an enumeration of }\psof{\omega}}$.
\end{xitemize}

\assert{4}: Suppose, toward a contradiction, that $\kappa\leq2^{\aleph_0}$. 
Let $\lambda>2^{\aleph_0}$ and let $\poP\in\calP$ be \st, for 
a $(\uniV,\poP)$-generic $\genG$, there are $j$, $M\subseteq\uniV[\genG]$ \st\
$\Elembed{j}{\uniV}{M}{\kappa}$, and $j(\kappa)>\lambda$. We have 
$M\models 2^{\aleph_0}\geq j(\kappa)$  by elementarity. Thus
$(2^{\aleph_0})^M\geq j(\kappa)>\lambda>(2^{\aleph_0})^\uniV$. Since
$(2^{\aleph_0})^\uniV=(2^{\aleph_0})^M$ by $\poP\in\calP$, this is a contradiction. 
\qedofLemma\qedskip

\begin{Thm}{\rm(Theorem 5.8 in \cite{sfetal-II})}\Label{p-MP-2-0}  Suppose that each element 
    of an iterable class $\calP$ of \pos\ is $\mu$-cc for some $\mu<\kappa$ and $\calP$ contains a \po\ $\poP$ which 
    adds a real. If $\kappa$ is tightly $\calP$-Laver-gen.\  
    superhuge then $\kappa=\continuum$.
\end{Thm}
\prf Suppose that $\kappa$ is tightly $\calP$-Laver-gen.\  
superhuge for the class of \pos\ $\calP$ as above. Then 
$\kappa\leq\continuum$ by \Lemmaof{p-Lg-RA-1-2-0},\,\assertof{3}. 

To prove $\continuum\leq\kappa$, let $\lambda\geq\kappa$, $2^{\aleph_0}$ be 
large enough.  By assumption there is a $\mu$-cc \po\ $\poQ$ \st\ there are $(\uniV,\poQ)$-generic
$\genH$ and $\elembed{j}{\uniV}{M}\subseteq\uniV[\genH]$ with\quad \wassertof{a} $\crit(j)=\kappa$,\\
\wassertof{b} $\cardof{\poQ}\leq j(\kappa)>\lambda$, \quad\wassertof{c} $\genH\in M$ and
\quad \wassertof{d} $j\imageof j(\kappa)\in M$. 

Since $\kappa$ is regular (\Lemmaof{p-Lg-RA-1-1}, \assertof{1}), and by elementarity, we 
have $M\modelof{j(\kappa)\mbox{ is regular}}$.  
By the closedness \wassertof{d} of $M$, it follows that $j(\kappa)$ is regular in
$\uniV[\genH]$. Hence it is also regular in $\uniV$. 

Thus, we have
$\uniV\modelof{j(\kappa)^{\mu}=j(\kappa)}$, since \SCH\ holds above $\kappa$ by 
Proposition 2.8,\,\assertof{1} in \cite{sfetal-II}. Since $\poQ$ has the $\mu$-cc and 
$\mu$, $\cardof{\poQ}\leq j(\kappa)$, it follows that
$\uniV[\genH]\modelof{\continuum\leq j(\kappa)}$. 
Again by \wassertof{d} (see \Lemmaof{L-lt-conti-0}\,\assertof{4}), 
we have $(j(\kappa)^+)^M=(j(\kappa)^+)^{\uniV}=(j(\kappa)^+)^{\uniV[\genH]}$. 
Thus $M\modelof{\continuum\leq j(\kappa)}$.

By elementarity, it follows that $\uniV\modelof{\continuum\leq\kappa}$.
\qedofThm\qedskip

\noindent
\prfof{\bfThmof{p-laver-1}} 
\assertof{A}: By \Lemmaof{p-Lg-RA-1-2-0}, \assertof{2} and \assertof{4}. \smallskip

\assertof{B}: By \Lemmaof{p-Lg-RA-1-2-0}, \assertof{2} and \assertof{3}. The last claim 
follows since $\MA(\calP)$ (and actually much more) holds by \Thmof{p-laver-2} below.

\assertof{$\Gamma$}: $\kappa$ is ``very large'' by \Lemmaof{p-Lg-RA-1-1}.
$\kappa\leq 2^{\aleph_0}$ follows from \Lemmaof{p-Lg-RA-1-2-0},\,\assertof{3}. The last 
statement follows from \Thmof{p-MP-2-0}.
\qedof{\Thmof{p-laver-1}}\qedskip\fi}

Laver-generic supercompactness also implies double plus versions of forcing axioms.
For a class $\calP$ of \pos\ and cardinals $\kappa$, $\mu$, let us denote with 
${\sf MA}^{+\mu}(\calP,\LT\kappa)$ and ${\sf MA}^{++\LT\mu}(\calP,\LT\kappa)$ the 
following versions of Martin's Axiom:

\begin{xitemize}
\item[\darkred$\darkred{\sf MA}^{+\mu}(\calP,\LT\kappa)$: ] 
  For any $\poP\in\calP$, 
  any family  
  $\calD$ of dense subsets of\/ $\poP$ with $\cardof{\calD}<\kappa$ and any family 
  $\calS$ of\/ $\poP$-names \st\ $\cardof{\calS}\leq\mu$ and
  $\forces{\poP}{\utilde{S}\xmbox{ is a stationary subset of }\omega_1}$ for all
  $\utilde{S}\in\calS$, there is a $\calD$-generic filter $\genG$ over $\poP$ \st\
  $\utilde{S}[\genG]$ is a stationary subset of $\omega_1$ for all
  $\utilde{S}\in\calS$. 

\item[\darkred$\darkred{\sf MA}^{++\LT\mu}(\calP,\LT\kappa)$: ]  
  For any $\poP\in\calP$, any family 
  $\calD$ of dense subsets of\/ $\poP$ with $\cardof{\calD}<\kappa$ and any family 
  $\calS$ of\/ $\poP$-names \st\ $\cardof{\calS}<\mu$ and
  $\forces{\poP}{\utilde{S}\xmbox{ is a stationary subset of }\Pkl{\eta_{\scriptstyle\utilde{S}}}{\theta_{\utilde{S}}}}$ 
  for some $\omega<\eta_{\utilde{S}}\leq\theta_{\utilde{S}}<\mu$ with $\eta_{\utilde{S}}$ regular, for all
  $\utilde{S}\in\calS$, there is a $\calD$-generic filter $\genG$ over $\poP$ \st\
  $\utilde{S}[\genG]$ is stationary in
  $\Pkl{\eta_{\scriptstyle\utilde{S}}}{\theta_{\utilde{S}}}$ for all 
  $\utilde{S}\in\calS$.  
\end{xitemize}

Clearly $\MA^{++\LT\omega_2}(\calP,\LT\kappa)$ is 
equivalent to $\MA^{+\omega_1}(\calP,\LT\kappa)$. 

{\ifextended
  {\extendedcolor
\begin{Thm}
  \Label{p-laver-2} For an iterable class $\calP$ of \pos, if $\kappa>\aleph_1$ 
  is $\calP$-Laver gen.\ supercompact then $\MA^{++\mu}(\calP,\LT\kappa)$ holds for all $\mu<\kappa$.
\end{Thm}
\prf Suppose that $\calP$ is an iterable class of \pos, $\kappa>\aleph_1$ 
is $\calP$-Laver gen.\ supercompact, and $\mu<\kappa$. Let $\calD$ and $\calS$ be as in 
the definition of $\MA^{++\mu}(\calP,\LT\kappa)$. \Wolog, we may assume that the underlying set of 
$\poP$ is some cardinal $\lambda_0$ and elements of $\calS$ are 
nice $\poP$-names.

Let $\lambda>\lambda_0$ be sufficiently large, and let $\utpoQ$ be a $\poP$-name \st\ 
$\forces{\poP}{\utpoQ\in\calP}$ and, for a $(\uniV,\poP\ast\utpoQ)$-generic filter $\genH$, 
there are   
transitive $M\subseteq\uniV[\genH]$ and $\Elembed{j}{\uniV}{M}{\kappa}$ with
\begin{xitemize}
\xitem[laver-0]{{}} 
  $j(\kappa)
  >\lambda$,
\xitem[laver-0-0]{{}} $\poP$, $\genH\in M$ and 
\xitem[laver-1]{{}} $j\imageof\lambda\in M$.
\end{xitemize}



By the choice of $\lambda$, \xitemof{laver-1} and \Lemmaof{L-lt-conti-0},\,\assertof{5}, we 
have $\poP$, $\calD$, $\calS\in M$. Let $\genG=\genH\cap\poP$. Then
$\genG\in M$ by \xitemof{laver-0-0}. 
Thus $\genG$ witnesses 

\begin{xitemize}
\xitem[laver-1-3] 
  $M\modelof{
  \begin{array}[t]{@{}l}
    \mbox{there is a }\calD\mbox{-generic filter }G\mbox{ over }\poP
    \\
    \mbox{\st\ }\utilde{S}(G)\mbox{ is a stationary subset of }\Pkl{\eta_{\scriptstyle\utilde{S}}}{\theta_{\utilde{S}}}
    \mbox{ for all }\utilde{S}\in\calS}.
  \end{array}$
\end{xitemize}
\iffalse{
This follows from \LemmaAof{laver-A-1} in case of \assertof{1} or 
from \xitemAof{laverA-20} in case of \assertof{2}.
}
\else 
\fi 

Since $j(\calD)=\setof{j(D)}{D\in\calD}$ and 
$j(\calS)=\setof{j(S)}{S\in\calS}$ by $\cardof{\calD}$, $\cardof{\calS}<\kappa$, 
$j(D)\supseteq j\imageof D$ for all $D\in\calD$, $j(S)\supseteq j\imageof S$ 
for all $S\in\calS$ and $j\imageof\genG\in M$ by 
\Lemmaof{L-lt-conti-0},\,\assertof{6}, 
it follows that 
\begin{xitemize}
\xitem[laver-1-4] 
  $M\modelof{
  \begin{array}[t]{@{}l}
    \mbox{there is a }j(\calD)\mbox{-generic filter }G\mbox{ over }j(\poP)
    \\
    \mbox{\st\ }\utilde{S}(G)\mbox{ is a stationary subset 
      of }\Pkl{\eta_{\scriptstyle \utilde{S}}}{\theta_{\utilde{S}}}
    \mbox{ for all }\utilde{S}\in j(\calS)}.
  \end{array}$
\end{xitemize}

By elementarity, it follows that 
\begin{xitemize}
\xitem[laver-1-5] 
  $\uniV\modelof{
  \begin{array}[t]{@{}l}
    \mbox{there is a }\calD\mbox{-generic filter }G\mbox{ over }\poP
    \\
    \mbox{\st\ }\utilde{S}(G)\mbox{ is a stationary subset of }\Pkl{\eta_{\scriptstyle\utilde{S}}}{\theta_{\utilde{S}}}
    \mbox{ for all }\utilde{S}\in \calS}.
  \end{array}$\vspace{-0.6ex}
\end{xitemize}\mbox{}\qedofThm
}
\else  
\begin{Thm}{\rm Theorem 5.7 in \cite{sfetal-II}}
  \Label{p-laver-2}
  \wassertof{1} For an iterable class $\calP$ whose elements are all 
  ccc, 
  if $\kappa$ 
  \memo{$\kappa>\aleph_1$ follows from this assumption.}
  is $\calP$-Laver-generically supercompact, then\quad
  $\MA^{++\LT\kappa}(\calP,\LT\kappa)$\quad holds.\smallskip
  
  \wassert{2} If\/ $\aleph_2$ is Laver-generically supercompact for an iterable  
  class $\calP$ of \pos, then\quad $\MA^{+\omega_1}(\calP)$\quad holds. \ifextended\else\qed\fi
\end{Thm}
{\ifextended\extendedcolor\prf \qedofThm\qedskip\fi}
\fi}

\begin{Prop}
  \Label{p-laver-3} If \ZFC\ $+$ ``there are two supercompact cardinals'' is consistent, 
  then \ZFC\ $+$ \FRP\ $+$ ``there is a tightly ccc-Laver-gen.\ supercompact'' is 
  consistent as well.
\end{Prop}
\prf Let $\kappa_0$ and $\kappa_1$ with $\kappa_0<\kappa_1$ be two supercompact cardinals.
We can use $\kappa_0$ to force $\MA^+(\sigma\mbox{-closed})$ by a \po\ of size $\kappa_0$. 
In the generic extension we have \FRP\ and $\kappa_2$ is still supercompact. Now we use 
$\kappa_1$ to force that $\kappa_1$ is tightly ccc-Laver-gen.\ supercompact in the generic 
extension as described in \Thmof{p-laver-0},\,\assertof{3}. \FRP\ still holds in the second 
generic extension since \FRP\ is preserved by ccc forcing (Theorem 3.4 in 
\cite{fjetal}).\qedofProp\qedskip

These results together with some other implications proved \cite{sfetal-II}
as well as some results that are going to be discussed \ifextended bellow \else in 
\cite{future} \fi
are integrated in \hyperref[figure1]{Figure 1} to obtain the following\ifextended\ extended 
diagram\fi:\newpage\bigskip

\hspace{-2em}
\includegraphics[width=1.0\textwidth]{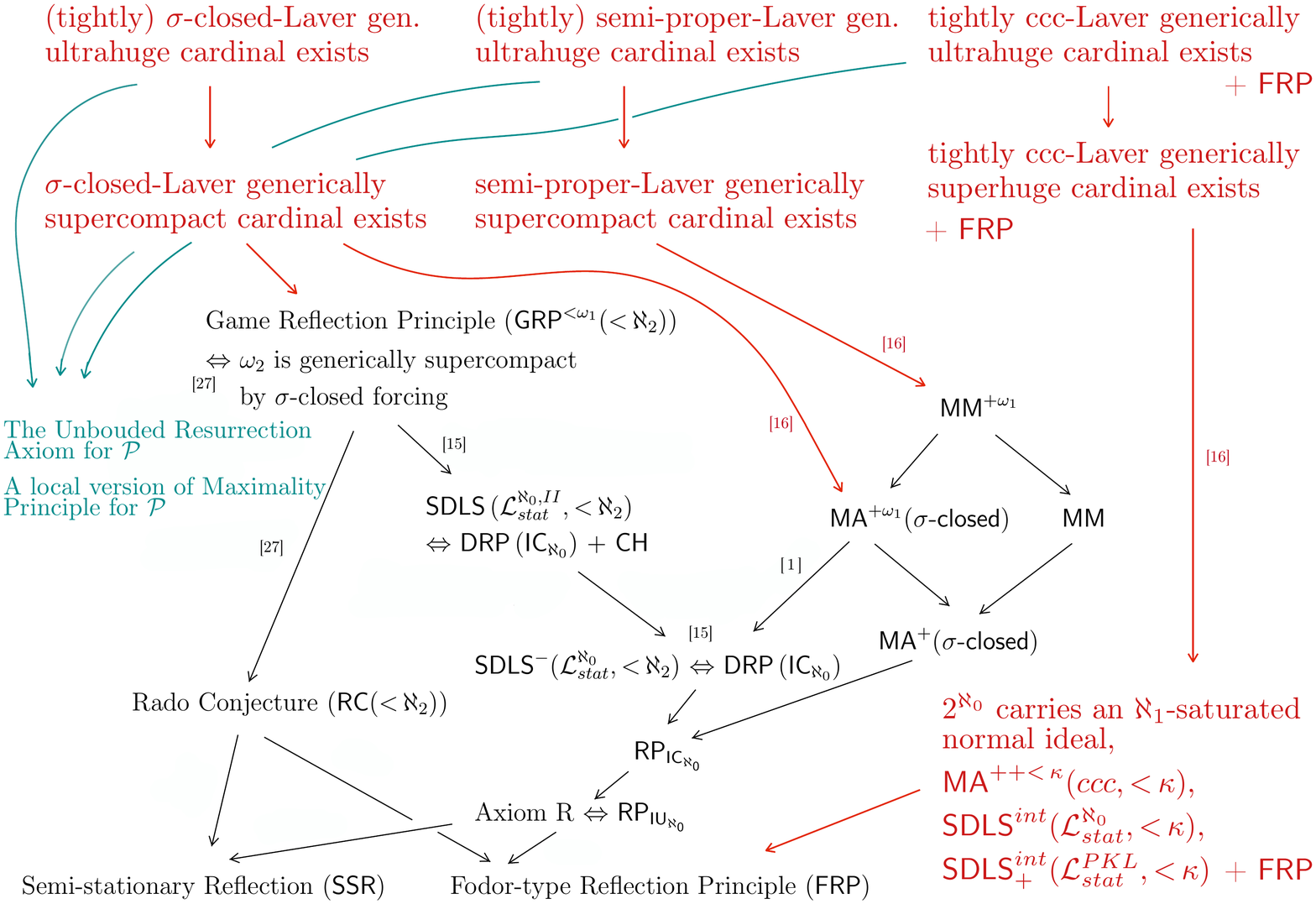}
\\[-2.52ex]
\mbox{}\hspace{19em}{\footnotesize\rotatebox[origin=c]{45}{$\darkred\Leftrightarrow$}}\\[-2.4\jot]
\mbox{}\hspace{7em}{\darkred\footnotesize many ``mathematical'' reflection theorems with reflection down to
  $\darkred\LT\aleph_2$ }\\[-2\jot]
\mbox{}\hspace{7em}{\darkred\tiny{[4]},
  [5], [8], [14], [22], etc.}
\vspace{-29.4ex}\\
{\tiny
\color[rgb]{0.0, 0.45, 0.55}
\Thmof{p-resurr-0}\\[-4\jot]
\Thmof{p-UR-4}}
\vspace{27.6ex}
\medskip\\
\phantomsection\Label{figure3}
\centerline{{\tt Figure 3.}}
\bigskip

For some iterable classes $\calP$ of \pos, even though they look quite natural, we can 
prove that there is no $\calP$-Laver generic large cardinal. 

\begin{Prop}
  \Label{p-Lg-RA-2} Suppose that
  $\calP$ is an iterable class of \pos\ \st\ all\/ $\poP\in\calP$ 
  are $\omega_1$-preserving and $\calP$ contains a \po\/ $\poP^*$ whose generic filter destroys a 
  stationary subset of $\omega_1$.\footnotemark Then there is no 
  $\calP$-Laver-gen.\ measurable cardinal.
\end{Prop}
\footnotetext{``$\poP^*$ destroys a stationary subset of $\omega_1$'' means here that a
  $\poP^*$-generic set codes a club subset of $\omega_1\setminus S$ in some absolute way. 
  \par 
  Note that, for stationary and co-stationary subset $S$ of $\omega_1$, various 
  posets are known which preserve $\omega_1$ while shooting a club in
  $\omega_1\setminus S$  (e.g.\ see \cite{abraham-shelah}).} 

\prf Suppose, toward a contradiction, that $\calP$ is as above and $\kappa$ is $\calP$-Laver 
gen.\ measurable cardinal.

Let $S\subseteq\omega_1$ be 
stationary (and co-stationary), and let $\poP^*\in\calP$ be a \po\ shooting a club in
$\omega_1\setminus S$. By assumption, there is 
a $\poP^*$-name $\utpoQ$ of a \po\ \st\ $\forces{\poP^*}{\utpoQ\in\calP}$ and, for
$(\uniV,\poP^*\ast\utpoQ)$-generic $\genH$, there are $j$, $M\subseteq\uniV[\genH]$ \st\ 
\begin{xitemize}
\xitem[x-Lg-RA-0] 
  $\Elembed{j}{\uniV}{M}{\kappa}$ and
\xitem[x-Lg-RA-1] $\poP$, $\genH\in M$.
\end{xitemize}
By the choice of $\poP^*$ ($\circleq\poP^*\ast\utpoQ$) and 
\xitemof{x-Lg-RA-1}, $M\modelof{S\mbox{ is a non-stationary subset of }\omega_1}$. Since
$\crit(j)=\kappa>\omega_1$ by \Lemmaof{p-Lg-RA-1-2-0},\,\assertof{1}\imemox{!!!}, 
we have $S=j(S)$. By $\uniV\modelof{S\xmbox{ is stationary subset of }\omega_1}$, this is a 
contradiction to the elementarity \xitemof{x-Lg-RA-0}.  
\qedofProp

\begin{Cor}
  Suppose that
  $\calP=\setof{\poP}{\poP\mbox{ is a \po\ preserving cardinals}}$ or
  $\calP=\setof{\poP}{\poP\mbox{ is an $\omega_1$-preserving \po\,}}$. Then there is no 
  $\calP$-Laver-gen.\ measurable cardinal.\qed
\end{Cor}

\begin{Prop}
  \Label{p-Lg-RA-3} Suppose that $\calP$ is an iterable class of \pos\ \st, for a regular 
  cardinal  $\delta>\aleph_1$, 
  \begin{xitemize}
  \xitem[x-Lg-RA-2] 
    all $\poP\in\calP$ are $\delta$-cc, and 
  \xitem[x-Lg-RA-3] 
    $\Fn(\omega,\omega_1,\LT\aleph_0)\in\calP$.\footnotemark 
  \end{xitemize}
  Then there is no $\calP$-Laver-gen.\ supercompact 
  cardinal.  
\end{Prop}
\footnotetext{Note that $\Fn(\omega,\omega_1,\LT\aleph_0)$ has $\omega_2$-cc (since its size 
  is $\aleph_1$). In particular It has the $\delta$-cc.}
\prf Suppose, toward a contradiction, that $\kappa$ is a $\calP$-Laver-gen.\ 
supercompact cardinal for $\calP$ as above. 

\begin{Claim}
  \Label{p-cl-Lg-RA-0} $\kappa=\omega_1$. 
\end{Claim}
\prfofClaim 
\memox{Trichotomy Theorem の証明を書き加えてそれのための Lemmata を書き加える．}
Suppose $\kappa\not=\omega_1$. Then, since $\kappa=\omega$ is impossible, we have
$\kappa>\omega_1$. 
Let $\poP:=\Fn(\omega,\omega_1,\LT\aleph_0)$. Since $\poP\in\calP$ by \xitemof{x-Lg-RA-3}, 
there is a 
$\poP$-name $\utpoQ$ of a \po\ \st, for $(\uniV,\poP\ast\utpoQ)$-generic $\genH$, there are $j$,
$M\subseteq\uniV[\genH]$ \st\ $\genH\in M$ and $\Elembed{j}{\uniV}{M}{\kappa}$. But then,  
by $\genH\in M$, we have $M\modelof{j(\omega_1)=\omega_1\mbox{ is countable}}$. This 
is a contradiction to the elementarity of $j$. \qedofClaim\qedskip

Let $\poP\in\calP$ be \st, for $(\uniV,\poP)$-generic $\genG$, there are $j$,
$M\subseteq\uniV[\genG]$ \st\ $\Elembed{j}{\uniV}{M}{\omega_1}$, and $j(\omega_1)>\delta$. 
Since $\poP$ is $\delta$-cc by \xitemof{x-Lg-RA-2},
$\uniV[\genG]\modelof{\delta\mbox{ is a cardinal}}$. It 
follows that 
$M\modelof{\delta\mbox{ is a cardinal and }\omega<\delta<j(\omega_1)}$. This is a 
contradiction to the elementarity of $j$. 
\qedofProp

\memox{Thanks to: Gappo, Sakai, Tsaprounis, Minden}

\section{Maximality Principle}
\Label{max} \ifextended{\extendedcolor{\It Maximality Principle} ({\It\MP}) in its non 
  parameterized version as given in  
Joel Hamkins' \cite{hamkins} was first formulated by Paul Larson following the ideas suggested  
by Christophe Chalons. 

In the language $\Lin$ of \ZFC, the Maximality Principle  
can only be formulated in 
an infinite set 
of $\Lin$-sentences asserting for each $\Lin$-sentence $\varphi$ that if it is a button then it is 
already pushed. }\else 
{\It Maximality Principle} ({\It\MP}) in its non parameterized form as introduced by 
Joel Hamkins in \cite{hamkins} is formulated in 
an infinite set 
of formulas asserting that all buttons are already pushed. 
\fi 
That is, for all $\Lin$-sentences
$\varphi$,\quad if, there is a \po\ $\poP$ \st\ 
\begin{xitemize}
\xitem[x-max-0] 
  $\forces{\poQ}{\varphi}$ holds for all \po\ $\poQ$ with $\poP\circleq\poQ$,\medskip
\end{xitemize}
then $\varphi$ holds.

If \xitemof{x-max-0} holds, then we shall say that $\varphi$ is a {\It button with the 
  push $\darkred\poP$}.{\ifextended\extendedcolor\footnote{\extendedcolor This paper was 
    written just when it was intensively discussed whether a dictator was going to push the 
    red button.}\fi}

 One of the easy consequences of \MP\ is the following:
\begin{Prop}{\rm(Hamkins \cite{hamkins})}
  \Label{p-max-0} \MP\ implies $\uniV\not=\uniL$. \qed
\end{Prop}
{\ifextended\extendedcolor Note that the statement ``$\uniV\not=\uniL$'' is apparently a button. 
For another consequence of \MP, see \Lemmaof{p-indep-0-4} below.
\fi}

For an  $\Lin$-sentence $\varphi$ let {\darkred $mp_\varphi$} be the $\Lin$-sentence:
\begin{xitemize}
\xitem[x-max-1] $\exists P\,(P\mbox{ is a \po}\land
  \forall Q(P\circleq Q\rightarrow\forces{Q}{\varphi}))\,\rightarrow\,\varphi$. 
\end{xitemize}

Formally we define \MP\ to be the collection of all $\Lin$-sentence of the 
form $mp_\varphi$ for $\Lin$-sentence $\varphi$. 

For an $\Lin$-sentence $\varphi$ let {\darkred$mp^+_\varphi$} be the $\Lin$-sentence:
\begin{xitemize}
\xitem[x-max-4] $\exists P\,(P\mbox{ is a \po}\land
  \forall Q(P\circleq Q\rightarrow\forces{Q}{\varphi}))\,\\ \rightarrow\,
  \forall R(R\mbox{ is a \po}\ \rightarrow \forces{R}{\varphi})$. 
\end{xitemize}
Let $\MP^+$ be the collection of $\Lin$-sentences of the form $mp^+_\varphi$ 
for all $\Lin$-sentences $\varphi$. 

\begin{Prop}\Label{p-max-4-1}{\rm( Hamkins \cite{hamkins})}
  \MP\ and $\MP^+$ are equivalent over $\ZFC$.
\end{Prop}
\prf It is clear that $\MP^+$ implies \MP.

To see that \MP\ implies $\MP^+$, let $\varphi$ be an arbitrary $\Lin$-sentence.
Let us write $\Box\varphi$ for
$\forall R(R\mbox{ is a \po}\ \rightarrow \forces{R}{\varphi})$.

It is easy to see that we have 
$\Box\varphi\ \leftrightarrow\ \Box\Box\varphi$. Thus $mp^+_\varphi$ is 
  equivalent to $mp_{\Box\varphi}$. The latter sentence is a member of \MP.
\qedofLemma\qedskip

{{\ifextended\extendedcolor\fi}
\begin{Lemma}
  \ifextended\extendedcolor\fi
  \Label{p-max-4-1-0} Suppose that $mp^+_\varphi$ holds for an $\Lin$-sentence $\varphi$. 
  Then for any \po\ $\poP$, \smash{$\forces{\poP}{mp^+_\varphi}$} holds.
\end{Lemma}
\prf This follows from the fact that the premise of $mp^+_\varphi$ is forcing absolute 
while the conclusion is forcing upward absolute. 
\qedofLemma\qedskip

\begin{Lemma}
  \ifextended\extendedcolor\fi
  \Label{p-max-1} Suppose that $\varphi_0$\ctentenc $\varphi_{n-1}$ are $\Lin$-sentences. 
  If \ZFC\ is consistent,  
  then so is \smash{\ZFC\ $+$ $mp^+_{\varphi_0}$ $+\cdots+$ $mp^+_{\varphi_{n-1}}$. }
\end{Lemma}
\prf Suppose otherwise. Then for 
some $\Lin$-sentences $\varphi_0$\ctentenc $\varphi_n$, \ZFC\ $+$ $mp^+_{\varphi_0}$
$+\cdots+$ $mp^+_{\varphi_n}$ is inconsistent. We can take $\varphi_0\ctentenc\varphi_n$ 
\st\ $n$ is minimal possible for such set of sentences. Then we have 
\begin{xitemize}
\xitem[p-max-4-2] $\ZFC$ $+$ $mp^+_{\varphi_0}$ $+\cdots+$ $mp^+_{\varphi_{n-1}}$ is 
  consistent and 

\xitem[x-max-2] 
  $\ZFC + mp^+_{\varphi_0} +\cdots+mp^+_{\varphi_{n-1}}  \vdash\neg mp^+_{\varphi_n}$. 
\end{xitemize}

Note that

\begin{xitemize}
\xitem[x-max-3] 
  $\neg mp^+_{\varphi_n}\ \leftrightarrow\
  \begin{array}[t]{@{}l}
    \exists P\,(P\mbox{ is a \po}\land
    \forall Q(P\circleq Q\rightarrow\forces{Q}{\varphi_n})) \\
    \land\ \exists R\,(R\mbox{ is a \po}\land\notforces{R}{\varphi_n}).
  \end{array}$
\end{xitemize}

In \ZFC\ $+$ $mp^+_{\varphi_0}$ $+\cdots+$ $mp^+_{\varphi_{n-1}}$, let $\poP$ be a \po\ 
whose existence is guaranteed by the first half of the right side of \xitemof{x-max-3}. 

Since all formulas of \ZFC\ $+$ $mp^+_{\varphi_0}$ $+\cdots+$ $mp^+_{\varphi_{n-1}}$ are 
forced by $\poP$ by \Lemmaof{p-max-4-1-0}, we have $\forces{\poP}{\neg mp^+_{\varphi_n}}$. 
Hence 
$\forces{\poP}{\exists R\,(\notforces{R}{\varphi_n})}$ by the second half of the right side 
of \xitemof{x-max-3}.

On the other hand, by the choice 
of $\poP$ we have $\forces{\poP}{\forall R\,(\forces{R}{\varphi_n})}$. 

Thus we have obtained a proof of contradiction from \ZFC $+$ $mp^+_{\varphi_0}$ $+\cdots+$
$mp^+_{\varphi_{n-1}}$. This is a contradiction to the 
assumption \xitemof{p-max-4-2}.
\qedofLemma\qedskip}






\begin{Thm}{\rm (Hamkins \cite{hamkins})}\Label{p-max-3}
 If \ZFC\ is consistent then so is \ZFC\ $+$ \MP.
\end{Thm}
\prf Assume that \ZFC\ $+$ \MP\ is inconsistent. Then there 
are $\Lin$-sentences $\varphi_0$\ctentenc\linebreak 
$\varphi_{n-1}$ \st\ \ZFC\ $+$ $mp^+_{\varphi_0}$ $+\cdots+$ $mp^+_{\varphi_{n-1}}$ is 
inconsistent (see \Propof{p-max-4-1}). This is a contradiction to \Lemmaof{p-max-1}.
\qedofThm\qedskip 

By practically the same argument as above, we can prove also the following:

\begin{Thm}
\Label{p-max-4}
  Suppose that ``x-large cardinal'' is a notion of a large cardinal formalizable in $\Lin$ \st, 
  \begin{xitemize}
  \xitem[x-max-4-0] 
    if $\kappa$ is 
    an x-large cardinal then the x-largeness of $\kappa$ is preserved by any set-forcing of 
    size $<\kappa$. 
  \end{xitemize}
  If \ZFC 
  $+$ ``there are class many x-large cardinals'' is consistent, then so is \ZFC\ $+$ \MP\
  $+$ ``there are class many x-large cardinals''.
\end{Thm}
\prf Working in the theory \ZFC\ $+$ ``there are class many x-large cardinals'' we have 
that $\forces{\poP}{\mbox{there are class many x-large cardinals}}$ holds for any \po\ $\poP$
since, by 
\xitemof{x-max-4-0}, only set many x-large cardinals are destroyed by $\poP$. 
Thus \Lemmaof{p-max-1} with \ZFC\ replaced by \ZFC\ $+$ ``there are class many x-large 
cardinals'' can be shown by the same argument. \qedofThm\qedskip

\begin{Thm}{\rm (Hamkins \cite{hamkins})}\Label{p-max-4-0}
  \MP\ is preserved by any set-generic extension.
\end{Thm}
{\ifextended\extendedcolor\prf By \Propof{p-max-4-1} and \Lemmaof{p-max-4-1-0}. \fi}
\qedofThm\qedskip

{\ifextended\extendedcolor
A sort of inverse of \Thmof{p-max-4} also holds: 

\begin{Thm}\extendedcolor\Label{p-max-5}{\rm( Hamkins \cite{hamkins})}
  Suppose that \MP\ holds. If ``x-large cardinal'' is a 
  notion of large cardinal formalizable in $\Lin$ \st
  \begin{xitemize}
  \xitem[x-max-5] 
    If $\kappa$ is an x-large cardinal, then $\kappa$ 
    is (weakly, resp.) inaccessible; and 
  \xitem[x-max-7] 
    no new x-large cardinal is 
    created by set-forcing. 
  \end{xitemize}
  If there is an x-large cardinal, then there are cofinally many x-large cardinals 
  in $\uniV$. 
\end{Thm}
\prf Suppose otherwise. Let $\kappa$ be an x-large cardinal, and
$\lambda>\kappa$ be a  
cardinal above which there are no x-large cardinals. 

Let $\poP$ be a \po\ which collapses $\lambda$ to be, say, of cardinality $\omega_1$, and let 
$\genG$ be a $(\uniV,\poP)$-generic filter. Then by \xitemof{x-max-5} 
and 
\xitemof{x-max-7}, 
there is no x-large cardinal in $\uniV[\genG]$. Also there is no x-large 
cardinal in 
any further generic extension by 
\xitemof{x-max-7}.

By \MP\ it follows that there is no x-large cardinal in $\uniV$. This is 
however a contradiction to the assumption of the theorem. 
\qedofThm\qedskip

If ``x-large cardinal'' implies (strong) inaccessibility, we can also prove \Thmof{p-max-5} 
by adding $\lambda$ many reals instead of collapsing cardinals below $\lambda$. This remark 
is going to be relevant in the proof of \Propof{p-para-max-1-0}.\imemox{!!!}.

In the following corollary, I call a cardinal $\kappa$ {\It resurrectably x-large} for a 
notion ``x-large'' of large cardinal, if there is a 
\po\ $\poP$ \st\ $\forces{\poP}{\kappa\mbox{ is x-large}}$.

\begin{Cor}\extendedcolor\Label{p-max-6}
  Suppose that \MP\ holds.\smallskip

  \wassert{1}
  If there is a (weakly, resp.) inaccessible cardinal then there are 
  class many (weakly, resp.) inaccessible cardinals.\smallskip

  \wassert{2} If ``x-large'' is a large cardinal property satisfying \xitemof{x-max-5} and 
  there is a resurresctably x-large cardinal then there are class many 
  resurrectably x-large cardinals. 
\end{Cor}
\prf \assertof{1}: The (weakly rep.) inaccessible cardinals satisfy  
\xitemof{x-max-5} and \xitemof{x-max-7}. \smallskip

\assertof{2} The notion of resurrectably x-large cardinal as a large cardinal property 
satisfies \xitemof{x-max-5} and \xitemof{x-max-7}.
\qedofCor
\fi}

\section{Independence of {\sf MP} under a Laver-gen.\ large cardinal}
\Label{indep}
In the following it is convenient to consider an abstract notion of large cardinal. As a 
generic name for a notion of large cardinal, we shall use the fancy words ``x-large 
cardinal'', ``y-large 
cardinal'' etc.\ which are in association with the German expression ``x-beliebig'' 
meaning ``really arbitrary''. This way of narration has been already used in the last section. 


Suppose that ``... is an x-large cardinal'' is a notion of large cardinal. 
We say this notion of large cardinal is {\It normal\/} if the following hold:
\begin{xitemize}
\xitem[x-max-7-0] ``$\kappa$ is an x-large cardinal'' is formalizable in $\Lin$ over \ZFC.
\xitem[x-max-8] 
  ``$\kappa$ is an x-large cardinal'' implies that $\kappa$ is 
  inaccessible;
\xitem[x-max-9] 
  ``$\kappa$ is an x-large cardinal'' cannot be destroyed by a forcing of size
  $\LT\kappa$;
\xitem[x-max-10] 
  No new x-large cardinal can be created by small forcing; and 
\xitem[x-max-11] 
  \ZFC\ + ``there are stationarily many x-large cardinals'' is consistent.\footnotemark
\end{xitemize}
\footnotetext{``there are stationarily many x-large cardinals'' is the axiom scheme 
  consisting of the statements ``if $C=\setof{\alpha\in\On}{\varphi(\alpha)}$ is a club in 
  $\On$ then there is an x-large cardinal $\kappa$ \st\ $\varphi(\kappa)$''
  for all $\Lin$-formulas $\varphi=\varphi(x)$. }

Note that most of the known notions of large cardinal are normal in the sense above under the 
assumption of the consistency of the existence of a sufficiently large cardinal.

\begin{Ex}\Label{ex-indep-a-0}
  The notion of ultrahuge cardinal is normal under the consistency of \ZFC\ $+$ ``there is a 
  2-almost-huge cardinal''.
\end{Ex}
\prf By Theorem 3.4 in Tsaprounis \cite{tsaprounis2}, if $\kappa$ is 2-almost-huge then 
there is a normal  
ultrafilter $\calU$ over $\kappa$ \st\
$\setof{\alpha<\kappa}{V_\kappa\modelof{\alpha\mbox{ is ultrahuge}}}\in\calU$. Thus 
$V_\kappa$ is a model of \ZFC $+$ ``there are stationarily many ultrahuge 
cardinals''.\qedofEx\qedskip 

\begin{Ex}
  \Label{ex-indep-a}
  The notion of super almost-huge cardinal is normal 
  under the consistency of \ZFC\ $+$ ``there is a huge cardinal''. 
\end{Ex}

The example above follows from the next theorem which should be a folklore:
\begin{Thm}\Label{p-indep-0}
  Suppose that $\kappa$ is huge. Then,
  $\setof{\alpha<\kappa}{V_\kappa\modelof{\alpha\xmbox{ is super almost-huge}}}$ is a 
  normal measure 1 subset of $\kappa$. \ifextended\else\qed\fi
\end{Thm}
{\ifextended\extendedcolor

Note that \Thmabove\ implies that $V_\kappa$ for a huge cardinal $\kappa$ models 
\ZFC\ $+$ there is a super almost-huge cardinal $+$ there are stationarily many 
inaccessible cardinals (actually there are stationarily many super almost-huge 
cardinals). See Theorems  
\ref{p-indep-0-6}, \ref{p-indep-0-7}.

\Thmof{p-indep-0} also tells that the existence of huge cardinal implies the consistency 
of the theory \ZFC\ $+$ $\delta$ is super almost-huge $+$ $V_\delta\prec\uniV$, c.f.\ 
\Thmof{p-para-max-0}.
\memox{!!!!}

We use the following \Thmof{p-indep-0-0} and \Lemmaof{p-indep-0-1} for the proof of 
\Thmof{p-indep-0}. 

For cardinals $\kappa\leq\lambda$ and a sequence
$\vec{\calU}=\seqof{\calU_\gamma}{\kappa\leq\gamma<\lambda}$ \st\  $\calU_\gamma$ is a 
normal ultrafilter over $\Pkl{}{\gamma}$ for all $\kappa\leq\gamma<\lambda$, we say
that $\vec{\calU}$ is {\darkred coherent} if 
$\calU_\gamma=\calU_\delta|\gamma
:=\setof{\setof{a\cap\gamma}{a\in A}}{A\in\calU_\delta}$ for all
$\kappa\leq\gamma\leq\delta<\lambda$. 

For a coherent sequence of normal ultrafilters
$\vec{\calU}=\seqof{\calU_\gamma}{\kappa\leq\gamma<\lambda}$, 
We let ${\darkred\Elembed{j_\gamma}{\uniV}{M_\gamma}{}}\cong Ult(\uniV,U_\gamma){\kappa}$ be 
the standard embedding, 
and, for $\kappa\leq\gamma\leq\delta<\lambda$, we define
$\darkred\Elembed{k_{\gamma,\delta}}{M_\gamma}{M_\delta}{}$ by
$k_{\gamma,\delta}([f]_{\calU_\gamma})
:=[\seqof{f(x\cap\gamma)}{x\in\Pkl{}{\delta}}]_{\calU_\delta}$.

Then we have $j_\delta=k_{\gamma\delta}\circ j_\gamma$. 

The following \Thmof{p-indep-0-0} is a slight modification of Theorem 24.11 in \cite{higher-inf}.
\begin{Thm}\Label{p-indep-0-0}\extendedcolor
  For a cardinal $\kappa$ and inaccessible $\lambda>\kappa$ \tfae:\smallskip

  \wassert{a} $\kappa$ 
  is an almost-huge cardinal with almost-huge elementary embedding $j$ with the target
  $j(\kappa)=\lambda$.  
  \smallskip

  \wassert{b} There is a coherent sequence $\seqof{\calU_\gamma}{\kappa\leq\gamma<\lambda}$ 
  of normal ultrafilters \st
  \begin{xitemize}
  \xitem[x-indep-a-0]
    for all $\kappa\leq\gamma<\lambda$ and $\alpha$ 
    with $\gamma\leq\alpha<j_\gamma(\kappa)$, there is $\gamma\leq\delta<\lambda$ \st\
    $k_{\gamma,\delta}(\alpha)=\delta$. \qed
  \end{xitemize}
\end{Thm}

\begin{Lemma}\extendedcolor
  \Label{p-indep-0-1}
  If $\kappa$ is an (almost) huge cardinal and
  \begin{xitemize}
  \xitem[x-indep-a-1] 
    $\Elembed{j}{\uniV}{M}{\kappa}$ is a(n almost) huge elementary embedding.
  \end{xitemize}

  Thus, in particular,
  \begin{xitemize}
  \xitem[x-indep-a-2] 
    $\fnsp{j(\kappa)\GT}{M}\subseteq M$. 
  \end{xitemize}
  Then\quad \wassertof{1} $j(\kappa)$ is inaccessible.\smallskip

  \wassert{2}
  $\setof{\alpha<\kappa}{\alpha\mbox{ is measurable}}$ is normal measure 1 subset 
  of $\kappa$.\smallskip
  
  \wassert{3} {$M\modelof{\setof{\alpha<j(\kappa)}{\alpha\mbox{ is measurable}}
      \mbox{ is stationary in }j(\kappa)}$.}\smallskip

  \wassert{4} $\setof{\alpha<j(\kappa)}{\alpha\mbox{ is measurable}}$ is cofinal 
  in $j(\kappa)$. 
\end{Lemma}
\prf
\imemox{/Users/Sakae\_1/TeX/talks/fuchino/kobe-set-theory-seminar-2023/kobe2023-06-05a.tex
  Lemma 2}
\assertof{1}: Since $\kappa$ is inaccessible.
$M\modelof{j(\kappa)\mbox{ is inaccessible}}$ by elementarity \xitemof{x-indep-a-1}.  By 
\xitemof{x-indep-a-2}, 
it follows that 
$j(\kappa)$ is really inaccessible.\smallskip

\assertof{2}:
$\kappa$ is measurable and an ultrafilter witnessing this 
is an element of $M$ by \xitemof{x-indep-a-2} and  \assertof{1}. Thus
$M\modelof{\kappa\mbox{ is measurable}}$. 
$\calU:=\setof{A\subseteq\kappa}{\kappa\in j(A)}$ 
is a normal ultrafilter over $\kappa$ and
$\setof{\alpha<\kappa}{\alpha\mbox{ is measurable}}\in\calU$.\smallskip

\assertof{3}: By \assertof{2}, $\setof{\alpha<\kappa}{\alpha\mbox{ is measurable}}$ is a 
stationary subset of $\kappa$. By elementarity \xitemof{x-indep-a-1}, it follows that
$M\modelof{\setof{\alpha< j(\kappa)}{\alpha{\mbox{ is measurable}}}
  \xmbox{ a is a stationary subset of } j(\kappa)}$.\smallskip

\assertof{4}: follows from \assertof{3} and 
\xitemof{x-indep-a-2}. \qedofLemma\qedskip

\noindent
\prfof{\bfThmof{p-indep-0}}
Let $\Elembed{j}{\uniV}{M}{\kappa}$ be a huge elementary embedding, so that we have
\begin{xitemize}
\xitem[x-indep-a-3] 
  $\fnsp{j(\kappa)}{M}\subseteq M$. 
\end{xitemize}

For $\kappa\leq\gamma<j(\kappa)$, let
$\calU_\gamma:=\setof{A\subseteq\Pkl{}{\gamma}}{j\imageof{\gamma}\in j(A)}$.
Then $\vec{\calU}:=\seqof{\calU_\gamma}{\kappa\leq\gamma<j(\kappa)}\in M$ by 
\xitemof{x-indep-a-3}, and 
$\vec{\calU}\models$\xitemof{x-indep-a-0}
(see the proof of 
  \cite{higher-inf}, Theorem 24.11).

Since 
\xitemof{x-indep-a-0}
is a closure property, $M$ knows that there are club many $\alpha<j(\kappa)$ 
\st\ $\seqof{\calU_\gamma}{\kappa\leq\gamma<\alpha}\models\mbox{
\xitemof{x-indep-a-0}}$ . 

By \Lemmaof{p-indep-0-1},\,(2), $M$ thinks that there are stationarily many 
$\alpha<\kappa$ which are inaccessible (actually even measurable!). Thus
{$M\modelof{\xmbox{there are stationarily many inaccessible }\alpha<j(\kappa)\xmbox{ \st\ }
  \seqof{\calU_\gamma}{\kappa\leq\gamma<\alpha}\models\mbox{
  \xitemof{x-indep-a-0}}}$}

By \Thmof{p-indep-0-0}, 
\begin{xitemize}
\xitem[x-indep-a-4] 
  $M\modelof{V_{j(\kappa)}\models\kappa\mbox{ is super almost-huge}}$.
\end{xitemize}

$\calU:=\setof{A\subseteq\kappa}{\kappa\in j(A)}$ is a normal ultrafilter over $\kappa$.
By 
\xitemof{x-indep-a-4} 
$\setof{\alpha<\kappa}{V_\kappa\modelof{\alpha\mbox{ is super almost-huge}}}\in\calU$. 
  \imemox{/Users/Sakae\_1/TeX/talks/fuchino/kobe-set-theory-seminar-2023/kobe2023-06-05a.tex
    Theorem 3}
\qedof{\Thmof{p-indep-0}}\qedskip

For an $\Lin$-formula $\psi=\psi(\overline{x})$, we shall call a large cardinal
$\kappa$ {\It$\psi$-absolute} if the formula $\psi$ is absolute between $V_\kappa$ and
$\uniV$ (i.e.\ if for any $\overline{a}\in V_\kappa$, we have
$V_\kappa\models\psi(\overline{a})$\ $\Leftrightarrow$\ $\uniV\models\psi(\overline{a})$, 
or more formally, if the $\Lin$-formula
$(\forall\overline{x}\in V_y)
(\psi^{V_y}(\overline{x}\leftrightarrow\psi(\overline{x})))$ holds for $y=\kappa$). 
\begin{Lemma}\extendedcolor
  \Label{p-indep-0-2} For any concretely given $n\in\natnums$, 
  there is an $\Lin$-formula $\psi^*_n$ \st\ for any inaccessible $\kappa$, $\kappa$ is
  $\psi^*_n$-absolute if and only if
  \begin{xitemize}
  \xitem[x-indep-a-5] for any $\calM\subseteq\uniV$ \st\ $\calM$ is a set forcing ground of
    $\uniV$ with $\uniV=\calM[\genG]$ where $\genG$ is an $(\calM,\poP)$-generic filter for 
    some \po\ $\poP\in (V_\kappa)^\calM$ (including the case of\/ $\poP=\ssetof{\bbone_\poP}$ and
    $\calM=\uniV$), we have that 
    all $\Sigma_n^\ZFC$-formulas are absolute between $(V_\kappa)^\calM$ and $\calM$.
  \end{xitemize}
\end{Lemma}
\prf By the analysis of set forcing ground in connection with Laver-Woodin theorem on 
definability of grounds (see e.g.\ \cite{reitz}). \qedofLemma\qedskip

\begin{Lemma}\extendedcolor
  \Label{p-indep-0-3} Suppose that $\psi^*_n$ is as in \Lemmaof{p-indep-0-2}. 
  Then $\psi^*_n$-absolute inaccessible cardinals are not resurrectable. I.e., if\/ $\poP$ is 
  a \pos\ and 
  $\forces{\poP}{\check{\lambda}\xmbox{ is }\psi^*_n\xmbox{-absolute inaccessible cardinal}}$, then 
  $\lambda$ is really $\psi^*_n$-absolute inaccessible cardinal. 
\end{Lemma}
\prf This is clear by the choice \xitemof{x-indep-a-5} of $\psi^*_n$. \qedofLemma\qedskip

\begin{Lemma}\extendedcolor
  \Label{p-indep-0-4} Assume that \MP\ holds. Suppose that $\psi^*_n$ for some 
$n\in\natnums$ is as in \Lemmaof{p-indep-0-2}. If there is a $\psi^*_n$-absolute 
  inaccessible cardinal, then there are unboundedly may $\psi^*_n$-absolute inaccessible cardinal.
\end{Lemma}
\prf 
``$\psi^*_n$-absolute inaccessible cardinal'' as an abstract notion of large cardinal 
satisfies \xitemof{x-max-5} and \xitemof{x-max-7}. Thus the Lemma follows from 
\Thmof{p-max-5}. 
\qedofLemma\qedskip

\begin{Lemma}\Label{p-indep-0-5}\extendedcolor
  Suppose that $\psi$ is an arbitrary $\Lin$-formula. If 
  there are stationarily many inaccessible cardinals, then there are cofinally 
  many $\psi$-absolute inaccessible cardinals.
\end{Lemma}
\prf Let $\lambda$ be an arbitrary cardinal. By Montague-Lévy Reflection Lemma 
\begin{xitemize}
\xitem[x-indep-a-6]
  $\calC:=\setof{\kappa\in\Card}{\kappa>\lambda,\,\psi\mbox{ is absolute between }
  V_\kappa\mbox{ and }\uniV}$
\end{xitemize}
contains a (definable) club subclass of $\Card$. By assumption there is an inaccessible $\kappa\in\calC$. Then
$\kappa>\lambda$ is $\psi$-absolute. \qedofLemma\qedskip

The following theorem says that there is no reasonable notion of large cardinal \st\ 
existence of that large cardinal implies \MP.
\begin{Thm}\Label{p-indep-0-6}\extendedcolor
  Suppose that ``x-large cardinal'' is a normal notion of large cardinal. 
  Then \ZFC\ $+$ ``\,there is an x-large cardinal'' $+$ $\neg\MP$ 
  is consistent. 
\end{Thm}
\prf Let $n\in\natnums$ be \st\ ``$\kappa$ is an x-large cardinal'' is $\Sigma^{\ZFC}_n$.

We work in \ZFC\ $+$ ``there is an x-large cardinals'' $+$ ``there are stationarily many 
inaccessible cardinals''. Note that this theory is consistent by the normality of the 
x-largeness.  

Let $\kappa$ be an x-large cardinal, and 
let $\kappa_0$ and $\kappa_1$ be the first two $\psi^*_n$-absolute inaccessible cardinals 
above $\kappa$ (they exist by \Lemmaof{p-indep-0-5}).

By $\psi^*_n$-absoluteness of $\kappa_1$ and the choice of $n$, we have
$V_{\kappa_1}\modelof{\kappa\xmbox{ is an x-large cardinal}}$. 
Since $\kappa_0$ is the unique $\psi^*_n$-absolute inaccessible cardinal in $V_{\kappa_1}$, 
we have $V_{\kappa_1}\not\models\MP$ by \Lemmaof{p-indep-0-4}. \qedofThm\qedskip

Similar theorem also holds for Laver-generic versions of normal notions of large cardinal.

\begin{Thm}\Label{p-indep-0-7}\extendedcolor
  Suppose that ``x-large cardinal'' is a normal notion of large 
  cardinal with Laver function and that its tight Laver-gen.\ version can be 
  forced similarly to \Thmof{p-laver-0} 
  for an iterable class $\calP$ of \pos\ given \Thmof{p-laver-0}. Then \ZFC\ $+$ ``there is a 
  tightly Laver gen.\ x-large cardinal for $\calP$'' $+$ $\neg\MP$ is consistent. 
\end{Thm}
\prf Let $n\in\natnums$ be \st\ the statement ``$\kappa$ is an x-large cardinal'' 
is $\Sigma^{\ZFC}_n$. 

Let $\kappa$ be 
an x-large cardinal, and let 
$\kappa_0$ 
and $\kappa_1$ be the first two $\psi^*_n$-absolute inaccessible cardinals 
above the $x$-large cardinal $\kappa$ (as before they exist by \Lemmaof{p-indep-0-5}). 

By the choice of $n$, $V_{\kappa_1}\modelof{\kappa\mbox{ is an x-large cardinal}}$. 
Thus by the assumption on the property ``x-large cardinal'', there is $\poP\in V_{\kappa_1}$ \st, for $(V_{\kappa_1},\poP)$-generic $\genG$, we have 
\begin{xitemize}
\xitem[] 
  $V_{\kappa_1}[\genG]\modelof{\kappa\mbox{ is tightly }\calP
  \mbox{-Laver-gen.\ x-large cardinal}}$.
\end{xitemize}

In $V_{\kappa_1}[\genG]$, $\kappa_0$ is still the unique $\psi^*_n$-absolute inaccessible 
cardinal above $\kappa$ by \Lemmaof{p-indep-0-3}. Thus $V_{\kappa_1}[\genG]\not\models\MP$ 
by \Lemmaof{p-indep-0-5}. 
\qedofThm

\fi}

{\ifextended\extendedcolor
\section{Boldface Maximality Principle for an iterable class $\calP$ of \pos\ and Laver-genericity}
\Label{para-max}
For an iterable class $\calP$ of \pos\ and (a definition of) a set $\Sigma$, the {\It Maximality 
Principle for $\calP$ with parameters from $\Sigma$} ($\darkred\MP(\calP,\Sigma)$) is the 
following principle: 

\begin{xitemize}
\item[{\darkred$\MP(\calP,\Sigma)$ }:] For 
  any $\Lin$-formula $\varphi=\varphi(x_0\ctenten)$ and $a_0\ctenten \in\Sigma$, if there 
  is $\poP\in\calP$ \st\ for any $\poP$-name $\utpoQ$ of a \po\ with 
  $\forces{\poP}{\utpoQ\in\calP}$, we have
  $\forces{\poP\ast\utpoQ}{\varphi(\check{a}_0\ctenten)}$, then we actually have
  $\forces{\poR}{\varphi(\check{a}_0\ctenten)}$ for all $\poR\in\calP$.\footnotemark
\end{xitemize}
\footnotetext{\extendedcolor In particular, since $\ssetof{\bbone}\in\calP$, 
  $\varphi(a_0\ctenten)$ holds. }

Similarly to \xitemof{x-max-0}, we shall call $\varphi(a_0\ctenten)$ as above a $\calP$-button, 
and $\poP$ a push of the $\calP$-button. 

\begin{Thm}\extendedcolor
  \Label{p-para-max-0} Suppose that ``x-large cardinal'' is a normal notion of large cardinal 
  with a Laver function \st\ the tight Laver-gen.\ version of x-large cardinal can be 
  forced similarly to \Thmof{p-laver-0} for one of the iterable classes $\calP$ of \pos\ 
  given in \Thmof{p-laver-0}.
  Working in $\uniV\modelof{\ZFC\ +\ \kappa\mbox{ is an x-large 
  cardinal }+\ V_\kappa\prec\uniV}$, there is a \po\ $\poP$ \st\ for 
  $(\uniV,\poP)$-generic $\genG$, we have 
  \begin{xitemize}
  \item[] 
    $\uniV[\genG]\modelof{\ZFC\ +\ \kappa\mbox{ is }\calP
    \mbox{-Laver gen.\ x-large cardinals }+\ \MP(\calP,\calH(\kappa))}$.
  \end{xitemize}
\end{Thm}

The following rather trivial lemma is used in the proof of \Thmabove.
\begin{Lemma}\extendedcolor
  \Label{p-para-max-1} If $\mapping{f}{\kappa}{V_\kappa}$ is a Laver function for an 
  x-large cardinal $\kappa$, then it is a book-keeping of elements of $V_\kappa$. I,e., for 
  any $a\in V_\kappa$ and $\alpha<\kappa$, there is $\beta\in\kappa\setminus\alpha$ \st\
  $f(\beta)=a$. 
\end{Lemma}
\prf Suppose that $a\in V_\kappa$ and $\alpha<\kappa$. Since $f$ is a Laver function for 
x-largeness of $\kappa$, there is $\Elembed{j}{\uniV}{M}{\kappa}$ with the closure property 
of $M$ corresponding to the x-largeness of $\kappa$ \st\ $j(f)(\kappa)=a=j(a)$. Note that 
we have $j(\alpha)=\alpha<\kappa<j(\kappa)$. 
By 
elementarity of $j$, it follows that there is $\beta\in\kappa\setminus\alpha$ \st\
$f(\beta)=a$.\qedofLemma\qedskip

\noindent
{\bf Proof of \bfThmof{p-para-max-0}:} We shall only consider the case that $\calP$ is the class of 
all proper posets. The other cases can be treated similarly. 

Let $\mapping{f}{\kappa}{V_\kappa}$ be a Laver function for x-largeness of $\kappa$.
Let
$\seqof{\poP_\alpha,\utildepoQ_\beta}{\alpha\leq\kappa,\beta<\kappa}$ be the  
CS-iteration of proper \pos\ defined by 
\begin{xitemize}
\xitem[x-para-max-0] 
  $\utildepoQ_\beta=\left\{
  \begin{array}{@{}l}
    f(\beta),\quad\ \ \mbox{if }f(\beta)\mbox{ is a }\poP_\beta\mbox{-name and }
    \forces{\poP_\beta}{f(\beta)\mbox{ is a proper \po}};\\
    \hfill\assertof{*}\\[\jot]
    \mbox{a }\poP_\beta\mbox{-name of a push of the }\calP\mbox{-button }\varphi(\uta_0\ctenten)
    \mbox{ in }V_\kappa,\\
    \qquad\qquad\mbox{if }f(\beta)\mbox{ is the }\Lin\mbox{-formula 
      with }\poP_\beta\mbox{-names }\uta_0\ctenten\in V_\kappa\mbox{, and}\\
    \qquad\qquad V_\kappa\modelof{\forces{\poP_\beta}{\varphi(\uta_0\ctenten)\mbox{ is 
          a }\calP\mbox{-button}}};\hfill\assertof{**}\\[\jot]
    \poP_\beta\mbox{-name of the trivial \po},\quad \mbox{otherwise}
  \end{array}
  \right.$
\end{xitemize}
for $\beta<\kappa$.

Let $\poP=\poP_\kappa$. 
Note that $\poP\in\calP$. Let $\genG$ be $(\uniV,\poP)$-generic. Then
$\uniV[\genG]\modelof{\kappa\mbox{ is }\calP\xmbox{-Laver gen. x-large}}$ by \assertof{*} 
in \xitemof{x-para-max-0} (see the proof of 
\Thmof{p-laver-0}). 


To show that $\uniV[\genG]$ satisfies $\MP(\calP,\calH(\kappa))$, 
let $a_0\ctenten\in\calH(\kappa)^{\uniV[\genG]}$ and $\Lin$-formula
$\varphi=\varphi(x_0\ctenten)$ be \st\  
\begin{xitemize}
\item[] 
  $\uniV[\genG]\models\varphi(a_0\ctenten)\mbox{ is a }\calP\mbox{-button}$.
\end{xitemize}
\Wolog, we may assume that $\uta_0\ctenten$ are $\poP$-names of $a_0\ctenten$ respectively 
and 
\begin{xitemize}
\xitem[x-para-max-1] 
  $\forces{\poP}{\varphi(\uta_0\ctenten)\mbox{ is a }\calP\mbox{-button}}$.\footnotemark
\end{xitemize}
\footnotetext{\extendedcolor E.g., by replacing $\varphi$ with 
  $((\varphi(x_0\ctenten)\land y\equiv0)\lor y\equiv1)$, and $\uta_0\ctenten$ with 
  $\uta_0\ctentenc\utb$. where $\utb$ is defined as follows: Let $A$ be a maximal antichain
  $\subseteq\setof{\condp\in\poP}{
    \condp\decides{\poP}{\varphi(\uta_0\ctenten)\xmbox{ is a button}} }$ (here, $\varphi$ is 
  the original $\varphi$ before the replacement), and 
  $\utb:=\setof{\pairof{\check0,\condp}}{\condp\in A,\,
    \condp\forces{\poP}{\varphi(\uta_0\ctenten)\xmbox{ is a button}}}\ \cup\ 
  \setof{\pairof{\check1,\condp}}{\condp\in A,\,
    \condp\forces{\poP}{\varphi(\uta_0\ctenten)\xmbox{ is not a button}}}$.}

Let $\alpha<\kappa$ be \st\ $\uta_0\ctenten$ are $\poP_\alpha$ names. 

For all $\alpha\leq\beta<\kappa$, since 
$\poP_\kappa\sim\poP_\beta\ast\utpoR$ where
$\forces{\poP_\beta}{\utpoR\mbox{ is proper}}$,  
\xitemof{x-para-max-1} implies that  we have
\begin{xitemize}
\xitem[x-para-max-2] 
   $\forces{\poP_\beta}{\varphi(\uta_0\ctenten)\mbox{ is a }\calP\mbox{-button}}$.
\end{xitemize}
By \Lemmaof{p-para-max-1}, there is $\alpha\leq\beta^*<\kappa$ \st\
$f(\beta^*)=\varphi(\uta_0\ctenten)$. By \xitemof{x-para-max-2}, and since
$V_\kappa\prec\uniV$, there is a push of the $\calP$-button $\varphi(\uta_0\ctenten)$ in
$V_\kappa$. Thus, by \assertof{$\ast\ast$} in \xitemof{x-para-max-0}, $\utpoQ_\beta$ is such a 
push. Since $\poP\sim\poP_{\beta^*}\ast\utpoQ_{\beta^*}*\utpoR$ where
$\forces{\poP_{\beta^*}\ast\utpoQ_{\beta^*}}{\utpoR\mbox{ is proper}}$, it follows that
$\forces{\poP}{\forall Q\,(Q\in\calP\ \rightarrow\forces{Q}{\varphi(\check{\uta}_0\ctenten)})}$.
Thus $\uniV[\genG]\models\forall Q\,(Q\in\calP\ \rightarrow\forces{Q}{\varphi(a_0\ctenten)})$. 
\imemox{Scan\_2023-06-16--18.56-resurrection - annotated.pdf\\!!!! pp.28--31}
\qedof{\Thmof{p-para-max-0}}\qedskip

The following proposition is a variation of \Thmof{p-max-5}

\begin{Prop}\Label{p-para-max-1-0}\extendedcolor Suppose that $\MP(\calP,\Sigma)$ holds for 
  an iterable class $\calP$ of \pos\ which contains either all \pos\ of the form
  $\Col(\omega_1,\lambda)$ or \pos\ adding arbitrary number of reals. Suppose further 
  that ``x-large cardinal'' is a notion of large cardinals formalizable in $\Lin$ \st\ 
  \begin{xitemize}
  \xitem[x-para-max-3] If $\kappa$ is an x-large cardinal, then $\kappa$ 
    is (weakly, resp.) inaccessible; and 
  \xitem[x-para-max-4] no new x-large cardinal is created by forcing by any $\poP\in\calP$. 
  \end{xitemize}
  If there is an x-large cardinal, then there are cofinally many x-large cardinals 
  in $\uniV$.
\end{Prop}
\prf Similarly to the proof of \Thmof{p-max-5}. See also the remark after the proof of
\Thmof{p-max-5}. \qedofProp\qedskip

For $n\in\natnums$, let $\psi^*_n$ be the $\Lin$-formula introduced in 
\Lemmaof{p-indep-0-2}. Since it is clear that ``$\psi^*_n$-absolute inaccessible cardinal'' 
a notion of large cardinal satisfying the conditions in \Propof{p-para-max-1-0}, we obtain 
the following:

\begin{Cor}\Label{p-para-max-1-1}\extendedcolor
  Suppose that $\MP(\calP,\Sigma)$ holds for 
  an iterable class $\calP$ of \pos\ which contains either all \pos\ of the form
  $\Col(\omega_1,\lambda)$ or \pos\ adding arbitrary number of reals. Then for any
  $n\in\natnums$, if there is a $\psi^*_n$-absolute inaccessible cardinal. Then there are 
  cofinally many $\psi^*_n$-absolute inaccessible cardinals. \qed
\end{Cor}

\begin{Thm}
  \Label{p-para-max-2}\extendedcolor
  Suppose that ``x-large cardinal'' is a normal notion of large cardinal with Laver function 
  such that the tight Laver-gen.\ version of x-large cardinal can be forced similarly to 
  \Thmof{p-laver-0} for one of the iterable classes $\calP$ of posets given in 
  \Thmof{p-laver-0}. 
  Working in
  $\uniV\modelof{\ZFC\ +\ \kappa\xmbox{ is an x-large cardinal }+\ V_\kappa\prec\uniV\ 
    +\xmbox{ there are }\psi^*_n\xmbox{-absolute inaccessible cardinals above }\kappa}$
  for sufficiently large $n$, and letting $\kappa_1$ the 
  second $\psi^*_n$-absolute inaccessible cardinal above $\kappa$, there is a poset $\poP$ 
  such that for $(V_{\kappa_1},\poP)$-generic $\genG$, we have 
  \begin{xitemize}
  \item[] 
    $\uniV_{\kappa_1}[\genG]\modelof{
    \begin{array}[t]{@{}l}
      \ZFC\ +\ \kappa\mbox{ is }
      \calP\mbox{-Laver gen.\,x-large cardinal}\\+\ \neg\MP\ +\ \neg\MP(\calP,\calH(\kappa))}.
    \end{array}$

  \end{xitemize}
\end{Thm}
\prf The proof of \Thmof{p-indep-0-6} works also here by \Corof{p-para-max-1-1}.
\qedofThm\qedskip

In spite of \Thmof{p-para-max-2}, the existence of a $\calP$-Laver generically ultrahuge 
cardinal implies a local version of maximality principle. 

To define the local version of maximality principle we are going to talk about below, let 
us call an $\Lin$-formula $\varphi=\varphi(x,a)$ with a parameter $a$ a {\It local 
  property  
  of cardinals} if, for any limit ordinal $\delta$ with $a\in V_\delta$ and a cardinal
$\mu<\delta$, we have $\big(V_\delta\models\varphi(\mu,a)\big)\ \leftrightarrow\ \varphi(\mu,a)$ and 
that this fact is provable in \ZFC\ ($+$ some formulas with the parameter $a$ which depict 
features of the set $a$). Being an inaccessible cardinal is a local property of 
cardinals, as well as being a Mahlo cardinal or being a measurable cardinal. In contrast,  
being a supercompact cardinal is not necessarily a local property of cardinals. 

A local property of cardinals $\varphi=\varphi(x,a)$ is a {\It local definition of a 
  cardinal} if there is provably at most one cardinal which satisfies the formula.

``The first inaccessible 
cardinal above a given cardinal $\mu$'' is a local definition of a cardinal as well as 
``the first measurable above $\mu$'' but not 
``the least supercompact cardinal above $\mu$''. 

If $\varphi(x,a)$ is a local definition of a cardinal, we denote the cardinal defined by
$\varphi(x,a)$ with {\It$\kappa^\bullet_{\varphi(x,a)}$}, {\It$\mu^\bullet_{\varphi(x,a)}$}, 
etc.\ or just with {\It$\kappa^\bullet$}, {\It$\mu^\bullet$}, etc.\ if 
we want to drop the explicit mention of the formula $\varphi(x,a)$ which defines the term. 
In the latter notation we identify the term $\kappa^\bullet$ with its definition
$\varphi(x,a)$ and say also that $\kappa^\bullet$ is a local definition of the cardinal. 

$\beth_\alpha(\omega_\beta)$ for any concretely 
given finite or countable ordinal $\alpha$, $\beta$ is another example of a local 
definition of a cardinal. 

Using this notation, we can show now that the existence of a $\calP$-Laver gen.\ ultrahuge 
cardinal implies the following local version of Maximality Principle for $\calP$. The 
following theorem is in line with the results in \cite{minden}:

\begin{Thm}\Label{p-UR-4}
  Suppose that $\calP$ is an iterable class of \pos\ and $\kappa$ is tightly $\calP$-Laver gen.\ 
  ultrahuge. Then, for any $\Lin$-formula $\varphi(x_0\ctentenc x_{n-1})$, $a_0$\ctentenc 
  $a_{n-1}\in\calH(\kappa)$, and a local definition $\mu^\bullet$ of a cardinal, if there is 
  $\poP\in\calP$ \st,
  \begin{xitemize}
  \xitem[x-UR-6] for any $\poP$-name $\utpoQ$ with $\forces{\poP}{\utpoQ\in\calP}$, we 
    have $\forces{\poP\ast\utpoQ}{V_{\mu^\bullet}
    \models\varphi(\check{a}_0\ctentenc\check{a}_{n-1})}$,
  \end{xitemize}
  then we have $(V_{\mu^\bullet})^\uniV\models\varphi(a_0\ctentenc a_{n-1})$. 
\end{Thm}
\prf Let $\kappa$, $\varphi$, $a_0$\ctentenc $a_{n-1}$, $\mu^\bullet$, $\poP$ as above.
Let $\lambda>(\mu^\bullet)^\uniV$ be a limit ordinal. Then there is a $\poP$-name $\utpoQ$ 
with $\forces{\poP}{\utpoQ\in\calP}$ \st, for $(V,\poP\ast\utpoQ)$-generic $\genH$, there 
are $j$, $M\subseteq\uniV[\genH]$ \st 
\begin{xitemize}
\xitem[x-UR-7] $\Elembed{j}{\uniV}{M}{\kappa}$,
\xitem[x-UR-8] $j(\kappa)>\lambda$, 
\xitem[x-UR-9] $\poP$, $\genH$, $(V_{j(\lambda)})^{\uniV[\genH]}\in M$, and 
\xitem[x-UR-10] $\poP\ast\utpoQ$ is forcing equivalent to a \po\ of size $j(\kappa)$. 
\end{xitemize}

By the choice of $\lambda$ and \xitemof{x-UR-7}, we have $j(\lambda)>(\mu^\bullet)^M$.
By \xitemof{x-UR-9} and \xitemof{x-UR-10}, we have $(V_{j(\lambda)})^M=(V_{j(\lambda)})^{\uniV[\genH]}$. 
Since $\mu^\bullet$ is a local definition, it follows that
$(\mu^\bullet)^M=(\mu^\bullet)^{\uniV[\genH]}$, and
$(V_{\mu^\bullet})^M=(V_{\mu^\bullet})^{\uniV[\genH]}$.  
Thus, by the choice of $\poP$, we have
$M\modelof{V_{\mu^\bullet}\models\varphi(a_0\ctentenc a_{n-1})}$. Since $a_i=j(a_i)$ for 
$i<n$ by \xitemof{x-UR-7}, it follows by the elementarity that
$(V_{\mu^\bullet})^\uniV\models\varphi(a_0\ctentenc a_{n-1})$. 
\qedofThm

\section{Resurrection Axioms}
\Label{resurr}
The following variants of Resurrection Axioms are introduced and studied by J.\ Hamkins 
and T.\ Johnstone (\cite{hamkins-johnstone}, \cite{hamkins-johnstone2}). 

For a class $\calP$ of \pos\ and a definition $\mu^\bullet$ of a cardinal (e.g.\ as 
  $\aleph_1$, $\aleph_2$, $\continuum$, $(\continuum)^+$. etc.) the {\It Resurrection Axiom for 
  $\calP$ and $\calH(\mu^\bullet)$} is defined by:
\begin{xitemize}
\item[$\darkred\RA^\calP_{\calH(\mu^\bullet)}$ :] For any
  $\poP\in\calP$, there is a $\poP$-name $\utpoQ$ 
  of \po\ \st\ $\forces{\poP}{\utpoQ\in\calP}$ and, for
  any $(\uniV,\poP\ast\utpoQ)$-generic $\genH$, we have
  $\calH(\mu^\bullet)^\uniV\prec\calH(\mu^\bullet)^{\uniV[\genH]}$. 
\end{xitemize}
 Here, $\mu^\bullet$'s in the left and right side of the last formula are actually meant
 $(\mu^\bullet)^\uniV$ and 
$(\mu^\bullet)^{\uniV[\genH]}$ respectively. 

The following boldface version of the Resurrection Axioms is also considered in 
\cite{hamkins-johnstone2}: 
For a class $\calP$ of \pos\ and a definition $\mu^\bullet$ of a cardinal (e.g.\ as 
  $\aleph_1$, $\aleph_2$, $\continuum$, $(\continuum)^+$. etc.) the {\It Resurrection Axiom 
  in Boldface for 
  $\calP$ and $\calH(\mu^\bullet)$} is defined by:

\begin{xitemize}
\item[$\darkred\BfRA^\calP_{\calH(\mu^\bullet)}$ :] For any $A\subseteq\calH(\mu^\bullet)$ and any 
  $\poP\in\calP$, there is a $\poP$-name $\utpoQ$ 
  of \po\ \st\ $\forces{\poP}{\utpoQ\in\calP}$ and, for
  any $(\uniV,\poP\ast\utpoQ)$-generic $\genH$, there is
  $A^*\subseteq\calH(\mu^\bullet)^{\uniV[\genH]}$ \st\ 
  $(\calH(\mu^\bullet)^\uniV,A,{\in})\prec(\calH(\mu^\bullet)^{\uniV[\genH]},A^*,{\in})$. 
\end{xitemize}
Clearly $\BfRA^\calP_{\calH(\mu^\bullet)}$ implies $\RA^\calP_{\calH(\mu^\bullet)}$. 

In the following we write ${\darkred\kappa_\refl}:=\max\ssetof{\aleph_2,2^{\aleph_0}}$. 
Note that this 
cardinal is the reflection point of the reflection properties we obtain in all scenarios of the 
trichotomy in \Thmof{p-laver-0} or \Thmof{p-laver-1}. 

\begin{Thm}\Label{T-Lg-RA-0}
  For an iterable class of \pos\ $\calP$, if
  $\kappa_\refl$ is tightly $\calP$-Laver-gen.\ superhuge, then
  $\BfRA^\calP_{\calH(\kappa_\refl)}$ holds.  

\end{Thm}
\prf The following proof is based on the idea suggested by Gunter Fuchs 
during a talk I gave at the New York Set Theory Seminar on October 7, 2022
. \memox{The proof should be rewritten because of $\RA\rightarrow\BfRA$. }



Suppose $A\subseteq\calH(\kappa_\refl)$ and $\poP\in\calP$. By the tightly $\calP$-Laver-gen.\  superhugeness 
of $\kappa_\refl$, there is a $\poP$-name $\utpoQ$ \smallskip of a \po\ with
$\forces{\poP}{\utpoQ\in\calP}$ \st, for $(\uniV,\poP\ast\utpoQ)$-generic $\genH$, there 
are $j$, $M\subseteq\uniV[\genH]$ with 
\begin{xitemize}
\xitem[was-a] $\Elembed{j}{\uniV}{M}{\kappa_\refl}$, 
\xitem[was-b] $j(\kappa_\refl)=\cardof{\poP\ast\utpoQ}$,
\xitem[was-c] $\poP$, $\genH\in M$, and
\xitem[was-d] $j\imageof{j(\kappa_\refl)}\in M$.
\end{xitemize}

\Wolog, we may assume that the underlying set of $\poP\ast\utpoQ$ 
is $j(\kappa_\refl)$. 

Since $\crit(j)=\kappa_\refl$, $j(a)=a$ for all $a\in(\calH(\kappa_\refl))^\uniV$. 

\begin{Claim}\Label{Cl-Lg-RA-0}
  $\calH(j(\kappa_\refl))^{\uniV[\genH]}\subseteq M$ and hence
  $\calH(j(\kappa_\refl))^M=\calH(j(\kappa_\refl))^{\uniV[\genH]}$\ .
\end{Claim}
\prfofClaim Suppose that $b\in\calH(j(\kappa_\refl))^{\uniV[\genH]}$ and let
$c\subseteq j(\kappa_\refl)$ be a code of $b$. Let $\utilde{c}$ be a nice 
$\poP\ast\utpoQ$-name of  $c$. By \xitemof{was-b}, $\cardof{\utilde{c}}\leq 
j(\kappa_\refl)$. By \xitemof{was-d}, it follows
that $\utilde{c}\in M$ (see \Lemmaof{L-lt-conti-0},\,\assertof{5}). Thus $c\in M$ by 
\xitemof{was-c} , and hence 
$b\in M$. \qedofClaim\smallskip 
 
Thus, we have
\begin{xitemize}
\item[] 
  $\Elembed{id_{\calH(\kappa_\refl)^\uniV}=j\restr\calH(\kappa_\refl)^\uniV}
  {\ (\calH(\kappa_\refl)^\uniV,A,{\in})\ }{\ (\calH(j(\kappa_\refl))^{\uniV[\genH]},j(A),{\in})}{}$.
\end{xitemize}
\qedofThm
\qedskip

The following strengthening of the Resurrection Axiom is introduced by Tsaprounis 
\cite{tsaprounis1}: 

For an iterable class $\calP$ of \pos, the {\It Unbounded Resurrection Axiom 
for $\calP$} is the following assertion. 

\begin{xitemize}
\item[{\darkred$\darkred\UR(\calP)$ : }] 
    For any $\lambda>\kappa_\refl$, and $\poP\in\calP$, there exists a $\poP$-name $\utpoQ$ 
    with $\forces{\poP}{\utpoQ\in\calP}$ \st, for $(\uniV,\poP\ast\utpoQ)$-gen.\ $\genH$, 
    there are $\lambda^*\in\On$ and $j_0\in\uniV[\genH]$ \st\ 
    $\Elembed{j_0}{\calH(\lambda)^\uniV}{\calH(\lambda^*)^{\uniV[\genH]}}{\kappa_\refl}$, and
    $j_0(\kappa_\refl)>\lambda$. 
\end{xitemize}
The ``tight'' version of the Unbounded Resurrection Axiom for
  $\calP$ will be also considered. 

\begin{xitemize}
\item[{\darkred$\TUR(\calP)$ : }] 
  For any $\lambda>\kappa_\refl$, and $\poP\in\calP$, there exists a $\poP$-name 
  $\utpoQ$ 
  with $\forces{\poP}{\utpoQ\in\calP}$ \st, for $(\uniV,\poP\ast\utpoQ)$-gen.\ $\genH$, 
  there are $\lambda^*\in\On$, and $j_0\in\uniV[\genH]$ \st\ 
  $\Elembed{j_0}{\calH(\lambda)^\uniV}{\calH(\lambda^*)^{\uniV[\genH]}}{\kappa_\refl}$, 
  $j_0(\kappa_\refl)>\lambda$, and {\darkblue$\darkblue\poP\ast\utpoQ$ is forcing equivalent to 
  a \po\ of size $j_0(\kappa_\refl)$.}
\end{xitemize}

Both of the principles can be yet extended to boldface versions similarly to the boldface 
version $\BfRA^\calP_{\calH(\mu^\bullet)}$ of $\RA^\poP_{\calH(\mu^\bullet)}$.  
However, $\UR(\calP)$ and 
$\TUR(\calP)$ can be easily proved to be equivalent to their respective boldface (apparent) 
extensions. 

\begin{Thm}\Label{p-resurr-0}
  For an iterable class $\calP$, if $\kappa_\refl$ is  
  tightly $\calP$-Laver gen.\ ultra\-huge, then 
    $\TUR(\calP)$  holds.
\end{Thm}
\prf Suppose that $\kappa_\refl$ is tightly $\calP$-Laver gen.\ ultrahuge.
Assume $\lambda>\kappa_\refl$, and $\poP\in\calP$. 

Let $\utpoQ$ be a $\poP$-name \st\ $\forces{\poP}{\utpoQ\in\calP}$ and, for
$(\uniV,\poP\ast\utpoQ)$-gen.\ filter $\genH$, there are $j,M\subseteq\uniV[\genH]$ \st\
\begin{xitemize}
\xitem[x-resurr-0] 
  $\Elembed{j}{\uniV}{M}{\kappa_\refl}$,
\xitem[x-resurr-1] 
  $j(\kappa_\refl)>\lambda$,
\xitem[x-resurr-2] 
  $\poP,\genH,V_{j(\lambda)}\in M$ and 
\xitem[x-resurr-3] $\poP\ast\utpoQ$ is forcing 
  equivalent to a \po\ of cardinality $j(\kappa_\refl)$. 
\end{xitemize}

\Wolog, let us assume that
\begin{xitemize}
\xitem[x-resurr-4] 
  $\cardof{\poP\ast\utpoQ}=j(\kappa_\refl)$, and $\poP\ast\utpoQ\subseteq V_{j(\kappa_\refl)}$. 
\end{xitemize}
Then $\calH(j(\lambda))^{\uniV[\genH]}\subseteq M$, since the code $\subseteq j(\lambda)$ of each 
element of $\calH(j(\lambda))^{\uniV[\genH]}$ has a $\poP\ast\utpoQ$-name in
$V_{j(\lambda)}\subseteq M$  
by \xitemof{x-resurr-3}, \xitemof{x-resurr-4} and \xitemof{x-resurr-2}. Thus the $\genH$ 
interpretation of 
the $\poP$-name of the code is in $M$ by \xitemof{x-resurr-2}. Hence the coded element 
of $\calH(j(\lambda))^{\uniV[\genH]}$ is also in $M$.  

It follows that $\calH(j(\lambda))^{M}=\calH(j(\lambda))^{\uniV[\genH]}$.

Thus, letting $j_0:=j\restr\calH(\lambda)^\uniV$, and $\lambda^*:=j(\lambda)$, we 
have
\begin{xitemize}
\item[] 
  $\Elembed{j_0}{\calH(\lambda)^\uniV}{
  \calH(\lambda^*)^{\uniV[\genH]}}{\kappa_\refl}$
  \quad and\quad $j_0(\kappa_\refl)=j(\kappa_\refl)>\lambda$ by \xitemof{x-resurr-1}.
\end{xitemize}
This shows that $\TUR(\calP)$ holds.
\qedofThm

\fi} 

\end{document}